\documentclass[a4paper,DIV=12,10pt,parskip=half]{scrartcl}

\usepackage[utf8]{inputenc}
\usepackage[USenglish]{babel}

\usepackage{hyperref}

\usepackage{amssymb,amsmath}

\usepackage[pdftex]{graphicx}

\usepackage{mathtools}
\mathtoolsset{centercolon}

\usepackage{multirow}
\usepackage{booktabs}

\usepackage{siunitx}

\usepackage{subcaption}
\usepackage{caption}

\usepackage{cleveref}

\usepackage{enumitem}

\usepackage{float}

\renewcommand{\vec}{\boldsymbol}
\newcommand{\mat}{\boldsymbol}

\renewcommand{\d}{\; \mathrm d}
\renewcommand{\P}{\mathbb{P}}
\newcommand{\R}{\mathbb{R}}
\newcommand{\N}{\mathbb{N}}

\newcommand{\force}{{\mat f}}

\newcommand{\QH}{Q_n^{\mathrm{H}}}
\newcommand{\IH}{I_{\tau}^{\mathrm{H}}}

\usepackage{amsthm}

\theoremstyle{definition}

\newtheorem{problem}{Problem}

\theoremstyle{remark}

\newtheorem{remark}{Remark}

\numberwithin{equation}{section}
\numberwithin{figure}{section}
\numberwithin{table}{section}
\numberwithin{problem}{section}

\begin{document}

\title{Higher order Galerkin--collocation time discretization with Nitsche's 
method for the Navier--Stokes equations}

\author{Mathias Anselmann$^{\star}$, Markus Bause$^{\star}$\\
	{\small ${}^\star$ Helmut Schmidt University, Faculty of 
		Mechanical Engineering, Holstenhofweg 85,}\\ 
	{\small 22043 Hamburg, Germany}\\
}
\date{}

\maketitle

\begin{abstract}
We propose and study numerically the implicit approximation in time of the Navier--Stokes equations by a Galerkin--collocation method in time combined with inf-sup stable finite element methods in space. The conceptual basis of the Galerkin--collocation approach is the establishment of a direct connection between the Galerkin method and the classical collocation methods, with the perspective of achieving the accuracy of the former with reduced computational costs in terms of less complex algebraic systems of the latter. Regularity of higher order in time of the discrete solution is ensured further. As an additional ingredient, we employ Nitsche's method to impose all boundary conditions in weak form with the perspective that evolving domains become feasible in the future. We carefully compare the performance poroperties of the Galerkin-collocation approach with a standard continuous Galerkin--Petrov method using piecewise linear polynomials in time, that is algebraically equivalent to the popular Crank--Nicholson scheme. The condition number of the arising linear systems after Newton linearization as well as the reliable approximation of the drag and lift coefficient for laminar flow around a cylinder (DFG flow benchmark with $Re=100$; cf.\ \cite{TS96}) are investigated. The superiority of the Galerkin--collocation approach over the linear in time, continuous Galerkin--Petrov method is demonstrated therein. 
\end{abstract}


\section{Introduction}

In the past, space-time finite element methods with continuous and 
discontinuous discretizations of the time and space variables have been studied strongly 
for the numerical simulation of incompressibe flow, wave propagation, transport phenomena 
or even problems of multi-physics; cf., e.g., 
\cite{ABM17,AM17,AB19_1,AB19_2,HST12,HST13,HST14,J93,KM04,S15,SY18,SY19}. Appreciable 
advantage of variational space-time discretizations is that they offer the 
potential to naturally construct higher order methods. In practice, these methods provide 
accurate results by reasonable numerical costs and on computationally feasible grids. 
Further, variational space-time discretizations allow to utilize fully adaptive finite 
element techniques to change the magnitude of the space and time elements in order to 
increase accuracy and decrease numerical costs; cf.\ \cite{BR12,BSB19,BBK20}. Strong relations 
between variational time discretization, collocation and Runge--Kutta methods have been 
observed \cite{AMN11,MN06}. Nodal superconvergence properties of variational time 
discretizations have also been proved \cite{BKRS18}. 

Recently, a modification of the standard continuous Galerkin--Petrov method (cGP) for 
the time discretization was introduced for wave problems (cf.~\cite{AB19_1,AB19_2,BMW17}). The modification comes through imposing collocation conditions involving the discrete solution's derivatives at the discrete time nodes while 
on the other hand reducing the dimension of the test space of the discrete variational problem compared with the standard cGP approach. A further key ingredient is the application of a special Hermite-type quadrature formula, proposed in~\cite{JB09}, and interpolation operator for the right-hand side function. Thereby, static condensation of degrees of freedom becomes feasible. These principles offer the potential of achieving the accuracy of Galerkin--Petrov methods with reduced computational costs by less complex algebraic systems. We refer to these schemes as Galerkin--collocation methods. Here, a family of Galerkin--collocation scheme with discrete solutions that are continuously differentiable in time and referred to as GCC$^1$ schemes is studied only. For families of schemes with even higher regularity in time we refer to \cite{AB19_2,BMW17}. The Galerkin--collocation schemes rely in an essential way on the perfectly matching set of the polynomial spaces (trial and test space), quadrature formula and interpolation operator. For wave problems, the GCC$^1$ approach has demonstrated its superiority over pure continuous Galerkin--Petrov approximations in time (cGP); cf.\ \cite{AB19_1}. Therefore, it seems to be natural to study the GCC$^1$ scheme also for the approximation of the Navier--Stokes equations. To the best of our knowledge, this paper is the first work, in that the application of such a Galerkin--collocation scheme to the Navier--Stokes equations is being proposed and investigated.

In the field of computational fluid dynamics, complex and dynamic geometries with moving boundaries are considered often. Fluid-structure interaction is a prominent example of multi-physics for problems with moving boundaries and interfaces. For these problems, Nitsche's fictitious domain method along with cut finite element techniques has been studied strongly in the recent past; cf.~\cite{BB16,BF07,BF19,BH14,BUM12,CH13,CRZ19,MLLR14,MSW18,S17} and the references therein. In this approach, the geometry is immersed into an underlying computational grid, which is not fitted to the geometric problem structure and, usually kept fixed over the whole simulation time. Thereby, mesh degeneration and remeshing are avoided for evolving and time-dependent domains. In Nitsche's method, the Dirichlet boundary conditions are added in a weak form to the variational equation of the partial differential equation, instead of imposing them on the definition of the function spaces of the variational problem. 

In this work, the Galerkin--collocation approximation of the Navier--Stokes equations is developed along with Nitsche's method for imposing Dirichlet boundary conditions. A Newton iteration is applied for solving the resulting nonlinear algebraic system of equations. For the Galerkin--collocation scheme GCC$^1(3)$ with piecewise cubic polynomials the algebraic system along with its Newton linearization is derived explicitly which is done here since the Galerkin--collocation approximation of the Navier--Stokes equations is presented for the first time and to facilitate the traceability of its implementation. The expected convergence behaviour of optimal order in time (and space) is demonstrated for the velocity and pressure variable. For flow around a cylinder it is illustrated that the accuracy of the approximation does not suffer from enforcing the boundary conditions in weak form by the application of Nitsche's method. Moreover, to show the superiority of the proposed approach over more standard time discretization schemes, a careful comparsion of the GCC$^1(3)$ approach with the continuous Galerkin--Petrov method using piecewise linear polynomials in time is performed. In algebraic form, the latter one can be recovered as the well-known Crank--Nicholson scheme. The errors of both approaches and, with regard to the future construction of efficient iterative solvers, the condition numbers of their Jacobian matrices are evaluated. Further, the performance properties for computing laminar flow around a cylinder are investigated. For this, the well-known DFG flow benchmark with $Re=100$ (cf.\ \cite{TS96}) as a challenging flow problem is used. The superiority of the GCC$^1(3)$ approach is cleary observed. Finally, we note that this work is considered as a building block for the future application of the proposed approach to fluid-structure interaction. In our outlook (cf.~Sec.~\ref{Sec:Outlook}), the feasibility for flow simulation on dynamically changing domains is demonstrated successfully. Here, a fictitious domain approach and stabilized cut finite element techniques are used; cf.~\cite{AB20}.

This paper is organized as follows. In Sec.~\ref{sec:Problem_Notation} we introduce our prototype model and notation. In Sec.~\ref{sec:GCC_with_Nitsche} we present its space-time discretization by Galerkin--collocation methods in time and inf-sup stable pairs of finite elements in space. Nitsche's method is applied to enforce Dirichlet 
boundary conditions in a weak sense. In Sec.~\ref{sec:C1_solution} the Galerkin--collocation scheme GCC$^1(3)$ of  piecewise cubic polynomials in time is studied. The algebraic formulation of the discrete system is derived and its solution by Newton's method is presented. In Sec.~\ref{Sec:NumExp} a careful numerical study of our approach 
is provided for the GCC$^1$(3) member of the Galerkin--collocation schemes.
In particular, the condition number of the arising algebraic systems after Newton linearization and the accuracy of the computed drag and lift coefficients for flow around a cylinder with $Re=100$ (DFG benachmark \cite{TS96}) are compared to the results of a standard cGP approximation with piecewise linear polynomials in time.
In Sec.~\ref{Sec:Outlook}, the potential of the presented approach, enriched by cut finite element techniques, for simulating flow on dynamically changing geometries by using fixed background meshes is illustrated. 

\section{Mathematical problem and notation}
\label{sec:Problem_Notation}

\subsection{Mathematical problem}

\begin{figure}[htb!]
	\centering
	\includegraphics[width=0.7\textwidth,keepaspectratio]
	{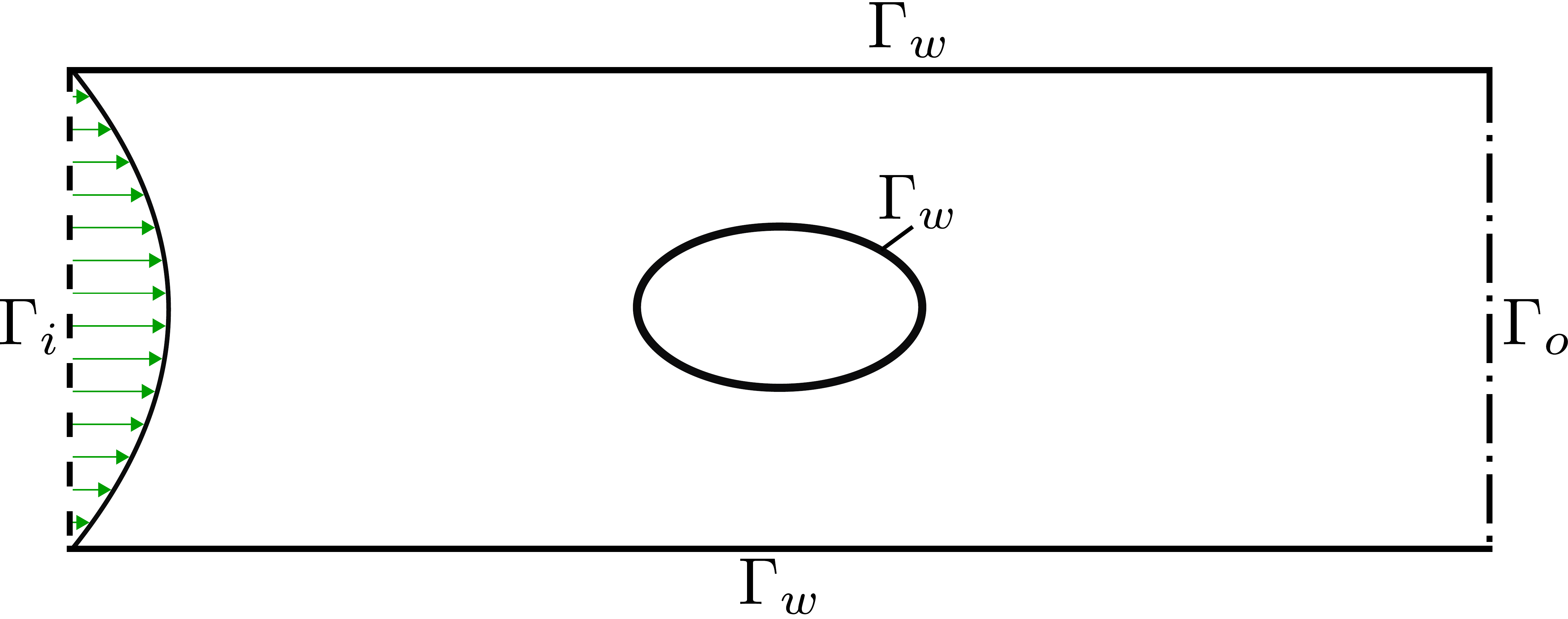}
	\caption{Notation for the flow domain $\Omega$.}
	\label{Fig:Domain}
\end{figure}

To fix our ideas and schemes in a familiar setting and simplify the notation, we 
restrict ourselves here to a common test problem, that of calculating nonstationary, 
incompressible flow past an obstacle, here taken as an inclined ellipse situated in a 
rectangle; cf.\ Fig.~\ref{Fig:Domain}. A generalization of our approach to more 
complex and three-dimensional bounded domains is straightforward. 

In this domain $\Omega \subset \R^2$ and for the time intervall $I=(0,T]$ we consider 
solving the Navier--Stokes equations, given in dimensionless form by 
\begin{subequations}
\label{eq:navier_stokes}
\begin{align}
\label{eq:navier_stokes_1}
\partial_t \vec{v} + (\vec{v} \cdot \nabla) \vec{v} - \nu \Delta \vec{v} + \nabla p	
& = \vec{f} &\quad& \hspace*{-3cm} \text{in } \Omega \times I\,,
\\
\label{eq:navier_stokes_2}
\nabla \cdot \vec{v} & = 0
&\quad& \hspace*{-3cm}\text{in } \Omega \times I\,,
\\
\label{eq:navier_stokes_3}
\vec{v} &= \vec{g} & & \hspace*{-3cm}\text{on } \Gamma_i \times I\,,
\\
\label{eq:navier_stokes_4}
\vec{v} &= \vec{0} & & \hspace*{-3cm}\text{on } \Gamma_w \times I\,,
\\
\label{eq:navier_stokes_5}
\nu \nabla \vec{v} \cdot \vec{n} - \vec{n}p &= \vec 0 & & \hspace*{-3cm}\text{on } 
\Gamma_o \times I\,,
\\
\label{eq:navier_stokes_6}
\vec{v}(0) &= \vec{v}_0 && \hspace*{-3cm}\text{in } \Omega\,.
\end{align}
\end{subequations}
In \eqref{eq:navier_stokes}, the unknows are the velocity field $\vec{v}$ and the 
pressure variable $p$. By $\nu$ we denote the dimensionless viscosity. Further, $\vec f$ 
is a given external force, $\vec v_0$ is the initial velocity and $\vec g$ is the 
prescribed velocity at the inflow boundary $\Gamma_i$. Eq.~\eqref{eq:navier_stokes_5} is 
the so-called do-nothing boundary condition for the outflow boundary $\Gamma_o$ with the 
outer unit normal vector $\vec n= \vec n(\vec x)$; cf.~\cite{J16}. At the upper 
and lower walls and on the boundary of the ellipse, jointly refered to as $\Gamma_w$, the 
no-slip boundary condition \eqref{eq:navier_stokes_4} is used.  For short, we put $\Gamma:= \Gamma_i\cup \Gamma_w\cup \Gamma_o$ and $\Gamma_D:=\Gamma_i\cup \Gamma_w$.

Wellposedness of the Navier--Stokes equations \eqref{eq:navier_stokes} and the existence 
of weak, strong or regular solutions to \eqref{eq:navier_stokes} in two and three space 
dimensions is not discussed here. The same applies to the optimal regularity of 
Navier--Stokes solutions at $t=0$ and the 
existence of non-local compatibility conditions. For the comprehensive discussion of these 
topics, we refer to the wide literature in this field; cf., e.g., \cite{J16,S01,T84} as well as 
\cite{B05,R83} and the references therein. Here, we assume the existence of a 
sufficiently regular (local) solution to the initial-boundary value problem 
\eqref{eq:navier_stokes} such that higher-order time and space discretizations become 
feasible. In particular we tacitly suppose that the solution to \eqref{eq:navier_stokes} is sufficient regular such that all of the equations given below are well-defined. 

\subsection{Notation}

In this work we use standard notation. $H^m(\Omega)$ is the Sobolev space of $L^2(\Omega)$ 
functions with derivatives up to order $m$ in $L^2(\Omega)$ and by $\langle \cdot,\cdot 
\rangle$ the inner product in $L^2(\Omega)$ and $(L^2(\Omega))^2$, respectively. In the 
notation of the inner product we do not differ between the scalar- and vector-valued case. 
Throughout, the meaning will be obvious from the context. We let 
\begin{equation*}
H^1_{0,\Gamma_D}(\Omega):=\{u\in  H^1(\Omega) \mid u=0 \mbox{ on } 
\Gamma_D=\Gamma_i\cup \Gamma_w \}\,.
\end{equation*}
For short, we put 
\begin{equation*}
Q:=L^2(\Omega)\,, \quad \vec V:= \left(H^1(\Omega)\right)^2\,,  \quad \vec V_0:= \left(H^1_{0,\Gamma_D}(\Omega)\right)^2 
\end{equation*}
and
\begin{equation*}
\vec X = \vec V \times Q\,, \quad  \vec X_0 := \vec V_0 \times Q \,.
\end{equation*}
Further, we define the function spaces 
\[
\vec V_{\text{div}} := \{\vec u\in \vec V\mid  \nabla \cdot \vec 
u = 0 \} \quad \text{and} \quad \vec V^0_{\text{div}} =  \vec V_{\text{div}}  \cap \vec V_0\,.
\]
We denote by $\vec V'$ the dual space of $\vec V_0$. 

For a Banach space $B$, we let $L^2(0,T;B)$, $C([0,T];B)$, and $C^m([0,T];B)$, $m\in\N$, be
the Bochner spaces of $B$-valued functions, equipped with their natural norms. 
For a subinterval $J\subseteq [0,T]$, we use the notations $L^2(J;B)$, $C^m(J;B)$, and 
$C^0(J;B):= C(J;B)$.

For the time discretization, we decompose the time interval $I=(0,T]$ into $N$ 
subintervals $I_n=(t_{n-1},t_n]$, $n=1,\ldots,N$, where $0=t_0<t_1< \cdots < t_{N-1} < 
t_N = T$ such that $I=\bigcup_{n=1}^N I_n$. We put $\tau = \max_{n=1,\ldots, N} \tau_n$ 
with $\tau_n = t_n-t_{n-1}$. Further, the set $\mathcal{M}_\tau := \{I_1,\ldots, I_N\}$
of time intervals is called the time mesh. For a Banach space $B$ and any $k\in \N_0$, 
we let 
\begin{equation*}
\P_k(I_n;B) = \Big\{w_\tau : I_n \to B \mid w_\tau(t) = \mbox{$\sum\limits_{j=0}^k$}
W^j t^j \; \forall t\in I_n\,, \; W^j \in B\; \forall j \Big\}\,.
\end{equation*}
For an integer $k\in \N$, we introduce the space
\begin{equation*}
X_\tau^k (B) := \left\{w_\tau \in C(\overline{I};B) \mid w_\tau|_{I_n} \in
\P_k(I_n;B)\; \forall I_n\in \mathcal{M}_\tau \right\}
\end{equation*}
of globally continuous functions in time and for an integer $l\in \N_0$ the space
\begin{equation*}
Y_\tau^{l} (B) := \left\{w_\tau \in L^2(I;B) \mid w_\tau|_{I_n} \in
\P_{l}(I_n;B)\; \forall I_n\in \mathcal{M}_\tau \right\}
\end{equation*}
of global $L^2$-functions in time.
For the space discretization, let $\mathcal{T}_h$ be a shape-regular mesh of $\Omega$ 
consisting of quadrilateral elements with mesh size $h>0$. For some  $r\in \N$, let 
$H_h=H_h^{(r)}$ be the finite element space given by 
\begin{equation}
\label{Eq:DefHh}
 H_h^{(r)}=\left\{v_h \in C(\overline{\Omega}) \mid v_h{}|_T\in \mathbb{Q}_r(K) \, 
\forall K \in \mathcal{T}_h \right\}\,,
\end{equation}
where $\mathbb{Q}_r(K)$ is the space defined by the multilinear reference mapping of polynomials on  the reference element with maximum degree $r$ in each variable. For brevity, we restrict our presentation to the Taylor--Hood family of inf-sup stable finite element pairs for the space discretization. The elements are used in the numerical experiments that are presented in Sec.\ \ref{Sec:NumExp}. However, these elements can be replaced by any other type of inf-sup stable elements. For some natural number $r\geq 2$ and with \eqref{Eq:DefHh} we then put 
\begin{equation*}
V_h = H_h^{(r)}\,, \quad V_h^0 = H_h^{(r)}\cap H^1_{0,\Gamma_D}(\Omega)\,, \quad  Q_h = H_h^{(r-1)}
\end{equation*}
and 
\begin{equation*}
X_h := V_h \times Q_h\,,\quad X_h^0 := V_h^0 \times Q_h\,,
\end{equation*}
as well as 
\begin{equation*}
\vec V_h = (V_h)^2 \,, \quad \vec V_h^0 = (V_h^0)^2\,, \quad \vec X_h = \vec V_h\times 
Q_h\,, \quad \vec X_h^0 = \vec V_h^0\times Q_h\,.
\end{equation*}
The space of weakly divergence free functions is denoted by
\begin{equation*}
 \vec V_h^{\text{div}} = \{\vec v_h \in \vec V_h \mid \langle \nabla \cdot \vec 
v_h,q_h\rangle  = 0 \; \text{for all } q_h \in Q_h\}\,.
\end{equation*}
For the discrete space-time functions spaces we use the abbreviations 
\begin{equation*}
\vec X^k_{\tau,h} = (X_\tau^k(V_h))^2\times  X_\tau^k (Q_h) \,, \quad \vec Y^k_{\tau,h} = (Y_\tau^k(V_h))^2\times  Y_\tau^k (Q_h) 
\end{equation*}
and 
\begin{equation*}
\vec X^{k,0}_{\tau,h} = (X_\tau^k(V_h^0))^2\times  X_\tau^k (Q_h) \,, \quad \vec Y^{k,0}_{\tau,h} = (Y_\tau^k(V_h^0))^2\times  Y_\tau^k (Q_h) 
\end{equation*}

Further, we define the semi-linear form $a: \vec X \times \vec X \rightarrow \R$ by 
\begin{equation}
\label{Eq:Slf_a}
\begin{aligned}
a(\vec u,\vec \phi) & :=  
\langle\partial_t \vec v, \vec \psi \rangle  + \langle (\vec v \cdot \vec \nabla) 
\vec v, \vec \psi \rangle + \nu \langle \nabla \vec v , \nabla \vec \psi  \rangle -\langle p, \nabla \cdot \vec \psi \rangle + \langle \vec \nabla 
\cdot \vec v, \xi \rangle 
\end{aligned}
\end{equation}
for $\vec u = (\vec v,p)\in \vec X$ and $\vec \phi = (\vec \psi,\xi) \in \vec X$ and the linear form $L:\vec X \rightarrow \R$ by 
\begin{equation*}
 L(\vec \phi; \vec f) := \langle \vec f, \vec \psi \rangle 
\end{equation*}
for $\vec \phi = (\vec \psi,\xi) \in \vec X$. For some parameters $\gamma_1>0$ and $\gamma_2> 0$ (to be discussed below) and 
\begin{equation*}
(\vec v\cdot \vec n)^- := \left\{\begin{array}{@{}ll} \vec v \cdot \vec n & \text{if}\;\; \vec v \cdot \vec n < 0\,,\\[1ex] 0 & \text{else}\,,\end{array}\right.
\end{equation*}
we introduce the semi-linear form $b: \vec H^{1/2}(\Gamma_D) \times 
\vec X_h \rightarrow \R$ by 
\begin{equation}
\label{Eq:Slf_bg}
\begin{aligned}
b_\gamma(\vec w,\vec \phi_h) : = &  - \langle \vec w, \nu \nabla \vec \psi_h \cdot \vec n + \xi_h  \vec n  \rangle_{\Gamma_D} - \langle ( \vec w\cdot \vec n)^-\vec w, \vec \psi_h \rangle_{\Gamma_D}  \\[1ex]
& \quad + \gamma_1 \nu \langle h^{-1} \vec w , \vec 
\psi_h  \rangle_{\Gamma_D}  + \gamma_2 \langle h^{-1} \vec w \cdot \vec n, \vec \psi_h \cdot \vec n 
\rangle_{\Gamma_D} 
\end{aligned}
\end{equation}
for $\vec w \in \vec H^{1/2}(\Gamma_D)$ and $\vec \phi_h = (\vec \psi_h,\xi_h) \in \vec X_h$. Finally, with \eqref{Eq:Slf_a} and \eqref{Eq:Slf_bg} the semi-linear form $a_\gamma:\vec X_h \times \vec X_h \rightarrow \R$ is given by 
\begin{equation}
\label{Eq:Slf_ag}
 a_\gamma (\vec u_h,\vec \phi_h) :=  a (\vec u_h,\vec \phi_h) - \langle \nu \nabla \vec 
v_h \cdot \vec n - p_h \vec n, \vec \psi_h\rangle_{\Gamma_D} + b_\gamma(\vec v_h,\vec 
\phi_h)
\end{equation}
for $\vec u_h = (\vec v_h,p_h)\in \vec X_h$ and $\vec \phi_h = (\vec \psi_h,\xi_h) \in 
\vec X_h$.

\section{Space-time finite element discretization with Ga\-lerkin--collocation time 
discretization and Nitsche's method}
\label{sec:GCC_with_Nitsche}

In this work, Nitsche's method \cite{N71} is applied within a space-time finite 
element discretization. In contrast to more standard formulations, the Dirichlet boundary 
conditions for the velocity field \eqref{eq:navier_stokes_3} and 
\eqref{eq:navier_stokes_4} are enforced weakly in the variational equation in terms 
of line integrals (surface integrals in three space dimensions). Our motivation for 
applying Nitsche's method comes through developing here a building block for flow 
problems with immersed or moving boundaries or even fluid-structure interaction that is based 
on using non-fitted background finite element meshes along with cut finite element 
techniques; cf.\ Sec.~\ref{Sec:Outlook}. For the time discretization a 
continuous Galerkin--Petrov approach (cf., e.g., \cite{HST12,HST13,HST14}) with discrete 
solutions $\vec v_{\tau,h}\in (C([0,T];X_h))^2$ and $p\in C([0,T];Q_h)$ is modified to a 
Galerkin--collocation approximation by combining the Galerkin techniques with the 
concepts of collocation. This approach has recently been developed \cite{AB19_1} and 
studied by an error analysis \cite{AB19_2} for wave equations. For this type of problems, 
the Galerkin--collocation approach has demonstrated its superiority over a pure Galerkin 
approach such that is seems to be worthwhile to apply the Galerkin--collocation technique 
also to the Navier--Stokes system \eqref{eq:navier_stokes}. 


\subsection{Space-time finite element discretization with Nitsche's fictious domain method}

A sufficiently regular solution of the Navier--Stokes system 
\eqref{eq:navier_stokes} satisfies the following variational space-time problem. 

\begin{problem}[Variational space-time problem]
\label{Prob:WS}
Let $\vec v_0 \in \vec V _{\operatorname{div}}$ be given. Let $\vec{\widehat g} \in \vec V_{\operatorname{div}}$ denote a prolongation of $\vec g$ such that $\vec{\widehat g}  = \vec g$ on $\Gamma_i$ and $\vec{\widehat g}  = \vec 0$ on $\Gamma_w$. Put $\vec{\widehat u} = (\vec{\widehat g},0)$. Let $\vec f \in L^2(0,T;\vec V')$ be given. Find $\vec u \in \vec{\widehat u}  +L^2 (0,T;\vec X_0)$ such that $\vec v(0) = \vec v_0$ and 
\begin{equation*}
 \int_0^T a(\vec u,\vec \phi) \d t  = \int_{0}^T  L(\vec \phi; \vec f) \d t
\end{equation*}
for all $\vec \phi \in L^2(0,T;\vec X_0)$
\end{problem}

For completeness and comparison, we briefly present the standard continuous 
Galerkin--Petrov approximation in time of Problem \ref{Prob:WS}, refered to as cGP($k$), 
along with the space discretization space; cf., e.g., \cite{HST12,HST13,HST14}. This 
reads as follows. 

\begin{problem}[Global problem of cGP($k$)]
	\label{Prob:DGS}
	Let an approximation $\vec v_{0,h} \in \vec V _h^{\operatorname{div}}$ of the initial value $\vec v_0$ be given. Let $\vec{\widehat g}_{\tau,h} \in \vec X^{k}_{\tau,h} $ denote a prolongation in the finite element spaces of the Dirichlet conditions on $\Gamma_i$ and  $\Gamma_w$ . Put $\vec{\widehat u}_{\tau,h} = (\vec{\widehat g_{\tau,h}},0)$. Let $\vec f \in L^2(0,T;\vec V')$ be given. Find $\vec u_{\tau,h} \in \vec{\widehat u} _{\tau,h} + \vec X^{k,0}_{\tau,h} $ such that $\vec v_{\tau,h}(0) = \vec v_{0,h}$ and 
	\begin{equation*}
	\int_0^T a(\vec u_{\tau,h},\vec \phi_{\tau,h}) \d t  = \int_{0}^T  L(\vec \phi_{\tau,h}; \vec f) \d t
	\end{equation*}
	for all $\vec \phi_{\tau,h} \in \vec Y^{k-1,0}_{\tau,h} $.
\end{problem}

By choosing test functions in $\vec Y^{k-1,0}_{\tau,h}$  supported on a single subinterval $I_n$ of the time mesh $\mathcal M_\tau$ we recast Problem~\ref{Prob:DGS}  as a time-marching scheme that is given by the following sequence of local problems on the subintervals $I_n$.

\begin{problem}[Local problem of cGP($k$)]
	\label{Prob:DLS}
	Let an approximation $\vec v_{0,h} \in \vec V _h^{\operatorname{div}}$ of the initial value $\vec v_0$ be given. Let $\vec{\widehat g}_{\tau,h} \in \vec X^{k}_{\tau,h} $ denote a prolongation into the finite element space of the Dirichlet conditions on $\Gamma_i$ and  $\Gamma_w$ . Put $\vec{\widehat u}_{\tau,h} = (\vec{\widehat g_{\tau,h}},0)$. Let $\vec f \in L^2(0,T;\vec V')$ be given. For $n=1,\ldots ,N$ and given $\vec u_{\tau,h}{}_{|I_{n-1}} \in \vec{\widehat u} _{\tau,h} + (\mathbb P_k(I_{n-1};V_h^0))^2 \times \mathbb P_k(I_{n-1};Q_h)$ find $\vec u_{\tau,h}{}_{|I_n} \in \vec{\widehat u} _{\tau,h} + (\mathbb P_k(I_n;V_h^0))^2 \times \mathbb P_k(I_n;Q_h)$ such that 
\begin{subequations}
\begin{align*}
\vec v_{\tau,h}{}_{|I_n} (t_{n-1})  & = \vec v_{\tau,h}{}_{|I_{n-1}} (t_{n-1})\,, \\[1ex]
p_{\tau,h}{}_{|I_n} (t_{n-1}) & = p_{\tau,h}{}_{|I_{n-1}} (t_{n-1})
\end{align*}
\end{subequations}
and 
	\begin{equation}
	\label{Eq:DLS:1}
	\int_{I_n} a(\vec u_{\tau,h},\vec \phi_{\tau,h}) \d t  = \int_{I_n} L(\vec \phi_{\tau,h}; \vec f) \d t
	\end{equation}
	for all $\vec \phi_{\tau,h} \in (\mathbb P_{k-1}(I_n;V_h^0))^2 \times \mathbb P_{k-1}(I_n;Q_h)$.
\end{problem}

In practice, the integral on the right-hand side of \eqref{Eq:DLS:1} is evaluated by 
means of an appropriate quadrature formula; cf.\ \cite{BKRS18,HST12,HST13,HST14}. 

\begin{remark}[Definition of initial pressure]
\label{Rem:IniPress}

\begin{itemize}
\hfill
\item In Problem ~\ref{Prob:DLS}, the quantities $\vec v_{\tau,h}{}_{|I_{n-1}}(t_{n-1})$ 
and $p_{\tau,h}{}_{|I_{n-1}}(t_{n-1})$ still need to be defined for $n=1$. For the 
velocity field we put $\vec v_{\tau,h}{}_{|I_{n-1}}(t_{n-1}):= \vec v_{0,h}$ for 
$n=1$ and with the approximation $\vec v_{0,h}$ of the initial value $\vec v_0$. Thus, it 
remains to define an approximation $p_{0,h}$ of the initial pressure $p_0:=p(0)$. 
This problem is more involved since the Navier--Stokes system does not provide an initial 
pressure. It is also impacted by the choice of the quadrature formula and the nodal 
interpolation properties of the temporal basis functions. A remedy based on Gauss 
quadrature in time and a post-processing for higher order pressure values in the 
discrete time nodes is proposed in \cite{HST12,HST14}. A further remedy consists in the 
application of a discontinuous Galerkin approximation (cf.\ \cite{HST14}) for the initial 
time step. In \cite{SR19}, a modification of the Crank--Nicholson scheme that is (up 
to quadrature) algebraically equivalent to the cGP(1) scheme is proposed by replacing the 
first two time steps with two implicit Euler steps. Regularity results for the Stokes 
equations, ensuring the optimal second order of convergence for the Crank--Nicholson 
scheme, are also studied in \cite{V94}. Nevertheless, this topic demands further research.

\item If the Navier--Stokes problem \eqref{eq:navier_stokes} is considered with  
(homogeneous) Dirichlet boundary conditions only, the unknown initial pressure 
$p_0:=p(0)\in L^2_0(\Omega)$ satisfies the boundary value probem (cf.~\cite[p.~376]{J16}, 
\cite{HR82})
\begin{subequations}
\label{Def:p0}
\begin{align}
\label{Def:p0_1}
- \Delta p_0 & = - \nabla \cdot \vec f(0) + \nabla \cdot ((\vec v_0 \cdot \nabla) \vec 
v_0) & & \hspace*{0cm} \text{in } \; \Omega\,, \\[1ex]
\label{Def:p0_2}
\nabla p_0 \cdot \vec n & = \big(\vec f(0) + \nu \Delta \vec v_0 \big)\cdot \vec n & & 
\hspace*{0cm}\text{on } \; \partial 
\Omega\,.
\end{align}
\end{subequations}	
In this case, we put $p_{\tau,h}{}_{|I_{n-1}}(t_{n-1}):= p_{0,h}$ for $n=1$ where 
$p_{0,h}\in Q_h\cap L^2_0(\Omega)$ denotes a finite element approximation of the solution 
$p(0)$ to \eqref{Def:p0}. 
\end{itemize}

\end{remark}

In the variational equation \eqref{Eq:DLS:1}, Dirichlet boundary conditions for the velocity field are enforced by the definition of the function space  $(\mathbb P_k(I_n;V_h^0))^2$ where by the discrete space $(V_h^0)^2$ homogeneous Dirichlet boundary conditions are prescribed for the velocity approximation $\vec v_{\tau,h}-\vec{\widehat g}_{\tau,h} $ and the test function $\vec \psi_{\tau,h}$. Using the Nitsche method, we solve instead of Problem \ref{Prob:DGS} the following one to that we refer to as cGP($k$)--N. 

\begin{problem}[Local Nitsche problem of cGP($k$): cGP($k$)--N]
		\label{Prob:DLSN}
		Let an approximation $\vec v_{0,h} \in \vec V _h^{\operatorname{div}}$ of the initial value $\vec v_0$ be given. For $n=1,\ldots ,N$ and given  $\vec u_{\tau,h}{}_{|I_{n-1}} \in (\mathbb P_k(I_{n-1};V_h))^2 \times \mathbb P_k(I_{n-1};Q_h)$  find $\vec u_{\tau,h}{}_{|I_n} \in (\mathbb P_k(I_n;V_h))^2 \times \mathbb P_k(I_n;Q_h)$ such that 
		\begin{equation}
		\label{Eq:DLSN:0}
		\begin{aligned}
		\vec v_{\tau,h}{}_{|I_n} (t_{n-1})  & = \vec v_{\tau,h}{}_{|I_{n-1}} (t_{n-1})\,, \\[1ex]
		p_{\tau,h}{}_{|I_n} (t_{n-1}) & = p_{\tau,h}{}_{|I_{n-1}} (t_{n-1})
		\end{aligned}
		\end{equation}
		and 
		\begin{equation}
		\label{Eq:DLSN:1}
		\int_{I_n} a_\gamma(\vec u_{\tau,h},\vec \phi_{\tau,h}) \d t  = \int_{I_n} 
 (L(\vec \phi_{\tau,h}; \vec f) + b_\gamma(\vec g,\vec \phi_{\tau,h})) \d t
		\end{equation}
		for all $\vec \phi_{\tau,h} \in (\mathbb P_{k-1}(I_n;V_h))^2 \times \mathbb P_{k-1}(I_n;Q_h)$.
\end{problem}

In Problem \ref{Eq:DLSN:1}, the Dirichlet boundary conditions for the velocity 
approximation are now ensured by the contribution of $b_\gamma$ to $a_\gamma$ and to the 
right-hand side of \eqref{Eq:DLSN:1}. For the definition of $a_\gamma$ and $b_\gamma$ we 
refer to \eqref{Eq:Slf_ag} and \eqref{Eq:Slf_bg}, respectively.  Let us still comment on 
the different boundary terms (line integrals in two dimensions and surface integrals in 
three dimensions) in the semi-linear forms \eqref{Eq:Slf_ag} and \eqref{Eq:Slf_bg}. The 
second term on the right-hand side of \eqref{Eq:Slf_ag} reflects the natural boundary 
condition, making the method consistent. The terms in $b_\gamma$ admit the following 
interpretation. The first term on the right-hand side of \eqref{Eq:Slf_bg} is introduced 
to preserve the symmetry properties of the continuous system. The second term 
incorporates the inflow condition. The last two term are penalty terms that insure the 
stability of the discrete system. In the inviscid limit $\nu = 0$, the last term amounts 
to a ''no-penetration'' condition. Thus, the semi-linear form $b_\gamma$ provides a 
natural weighting between boundary terms corresponding to viscous effects ($\vec v = \vec 
g$), convective behaviour ($(\vec v\cdot \vec n)^- \vec v =  (\vec g\cdot \vec n)^- \vec 
g$) and inviscid behaviour ($\vec v \cdot \vec n = \vec g \cdot \vec n$). Since the 
second term on the right-hand side of \eqref{Eq:Slf_bg} introduces a further 
nonlinearity, this term with little influence is ignored in the case of low Reynolds 
number flow that is assumed here (cf.~Sec.~\ref{Sec:NumExp}).
 
\subsection{Galerkin--collocation time discretization}
\label{sec:galerkin_collocation_time_discretization}

Our modification of the standard continuous Galerkin--Petrov method (cGP) for 
time discretization, that is used in Problem \ref{Prob:DLS}, and the innovation of this 
work comes through imposing a collocation condition involving the discrete solution's 
first derivative at the endpoint of the subinterval $I_n$ along with $C^1$-continuity 
constraints at the initial point of the subinterval $I_n$ while on the other hand 
downsizing the test space of the variational equation \eqref{Eq:DLS:1}. This principle is applied with the perspective of achieving the accuracy of Galerkin schemes with reduced computational costs. We refer to this family of schemes combining Galerkin and collocation techniques as Galerkin--collocation methods, for short GCC$^1$($k$), 
where $k$ denotes the degree of the piecewise polynomial approximation in time and the 
part C$^1$ in GCC$^1$($k$) denotes the continuous differentiability of the discrete 
solution. A further key ingredient in the construction of the Galerkin--collocation approach 
comes through the application of a special quadrature formula, investigated in~\cite{JB09}, and the 
definition of a related interpolation operator for the right-hand
side term of the variational equation. Both of them use derivatives of the given 
function. The Galerkin--collocation schemes rely in an essential way on the perfectly 
matching set of the polynomial spaces (trial and test space), quadrature formula, and 
interpolation operator. In particular, a condensation of degrees of freedom becomes feasibel by the construction principle such that smaller algebraic systems are obtained. The concept of Galerkin--collocation approximation was recently introduced in  \cite{BMW17} for systems of ordinary differential equations and applied successfully to 
wave problems in \cite{AB19_1,AB19_2}. Besides the numerical studies given in \cite{AB19_1}, showing the 
superiority of the Galerkin-collocation approach over a pure Galerkin approach as used in 
Problem \ref{Prob:DLS}, a rigorous error analysis is provided for the 
Galerkin--collocation approximation of wave phenomena in \cite{AB19_2}. 

From now on we assume a polynomial degree of $k\ge 3$. To introduce the 
Galerkin--collocation approximation, we need to define the Hermite quadrature formula and 
the corresponding interpolation operator. Let $\hat{t}^{\,\mathrm{H}}_1=-1$, 
$\hat{t}^{\,\mathrm{H}}_{k-1}=1$, and $\hat{t}^{\,\mathrm{H}}_s$,
$s=2,\dots,k-2$, be the roots of the Jacobi polynomial on $\widehat{I}:=[-1,1]$ with
degree $k-3$ associated to the weighting function $(1-\hat{t})^2(1+\hat{t})^2$.
Let $\widehat{I}^{\,\mathrm{H}}:C^1\big(\widehat{I};B\big)\to \P_k\big(\widehat{I};B\big)$
denote the Hermite interpolation operator with respect to point value and first
derivative at both $-1$ and $1$ as well as the point values at $\hat{t}_s^{\,\mathrm{H}}$,
$s=2,\dots,k-2$. By
\begin{equation*}
\widehat{Q}^{\mathrm{H}}(\hat{g}) := \int_{-1}^{1}
\widehat{I}^{\,\mathrm{H}}(\hat{g})(\hat{t})\d \hat{t}
\end{equation*}
we define an Hermite-type quadrature on $[-1,1]$ which can be written as
\begin{equation*}
\widehat{Q}^\mathrm{H}(\hat{g}) = \widehat{\omega}_L \hat{g}'(-1) + \sum_{s=1}^{k-1}
\widehat{\omega}_s \hat{g}(\hat{t}^{\,\mathrm{H}}_s) + \widehat{\omega}_R \hat{g}'(1)\,,
\end{equation*}
where all weights are non-zero. Using the affine mapping $T_n:\widehat{I}\to \overline{I}_n$
with $T_n(-1) = t_{n-1}$ and $T_n(1) = t_n$, we obtain
\begin{equation}
\label{Eq:GLHF}
\QH(g) = \left(\frac{\tau_n}{2}\right)^2 \widehat{\omega}_L \d_t g(t_{n-1}^+)
+ \frac{\tau_n}{2} \sum_{s=1}^{k-1} \widehat{\omega}_s g(t^\mathrm{H}_{n,s})
+ \left(\frac{\tau_n}{2}\right)^2 \widehat{\omega}_R \d_t g(t_n^-)
\end{equation}
as Hermite-type quadrature formula on $I_n$, where $t^\mathrm{H}_{n,s} :=
T_n(\hat{t}^\mathrm{H}_s)$, $s=1,\dots,k-1$. We note that $\QH$ as 
defined by \eqref{Eq:GLHF} integrates all polynomials up to degree $2k-3$ exactly, 
cf.~\cite{JB09}. Using $\widehat{I}^{\,\mathrm{H}}$ and $T_n$, the local Hermite 
interpolation on $I_n$ is given by
\begin{equation}
\label{Def:I_n}
I_n^\mathrm{H} : C^1(\overline{I}_n;B)\to \P_k(\overline{I}_n;B) \,,
\qquad 
v\mapsto \big(\widehat{I}^{\,\mathrm{H}}(v\circ T_n)\big)\circ T_n^{-1}\,.
\end{equation}
Moreover, for all $n=1,\dots,N$ we define the global Hermite
interpolation $\IH:C^1(\overline{I};B)\to X_\tau^k(B)$ by means of 
\begin{equation}
\label{Eq:DefHIntOp}
\IH w|_{I_n} := I_n^{\mathrm{H}}(w|_{I_n})\,.
\end{equation}
This operator is applied componentwise to vector-valued functions.

The local problem of the Galerkin--collocation approach along with Nitsche's method 
for enforcing Dirichlet boundary conditions then reads as follows.

\begin{problem}[Local $I_n$ problem of GCC$^1$($k$) with Nitsche's method: GCC$^1$($k$)--N]
\label{Prob:DLSNGC}
Let $k\geq 3$ and an approximation $\vec v_{0,h} \in \vec V _h^{\operatorname{div}}$ of 
the initial value $\vec v_0$ be given. For $n=1,\ldots ,N$ and  given $\vec u_{\tau,h}{}_{|I_{n-1}} 
\in (\mathbb P_k(I_{n-1};V_h))^2 \times \mathbb P_k(I_{n-1};Q_h)$ find $\vec u_{\tau,h}{}_{|I_{n}} 
\in (\mathbb P_k(I_n;V_h))^2 \times \mathbb P_k(I_n;Q_h)$ such that
\begin{align}
\label{DLSNGC:1_1}
\vec v_{\tau,h}{}_{|I_n} (t_{n-1})  & = \vec v_{\tau,h}{}_{|I_{n-1}} (t_{n-1})\,,
& \partial_t \vec v_{\tau,h}{}_{|I_n} (t_{n-1})  & = \partial_t 
\vec v_{\tau,h}{}_{|I_{n-1}} (t_{n-1})\,, \\[1ex]
\label{DLSNGC:1_2}
p_{\tau,h}{}_{|I_n} (t_{n-1}) & = p_{\tau,h}{}_{|I_{n-1}} (t_{n-1})\,,
& \partial_t p_{\tau,h}{}_{|I_n} (t_{n-1}) & = \partial_t p_{\tau,h}{}_{|I_{n-1}} 
(t_{n-1})\,,
\end{align}
and
\begin{equation}
\label{Eq:DLSNGC:2}
a_\gamma(\vec u_{\tau,h}(t_n),\vec \phi_{h}) = L(\vec \phi_{h}; \IH \vec f(t_n)) + 
b_\gamma(\IH \vec g(t_n),\vec \phi_{h})
\end{equation}
for all $\vec \phi_h \in \vec X_h$ as well as 	
\begin{equation}
\label{Eq:DLSNGC:3}
\int_{I_n} a_\gamma(\vec u_{\tau,h},\vec \phi_{\tau,h}) \d t  = \int_{I_n}  (L(\vec \phi_{\tau,h}; \IH\vec f) + b_\gamma(\IH \vec g,\vec \phi_{\tau,h})) \d t
\end{equation}
for all $\vec \phi_{\tau,h} \in (\mathbb P_{k-3}(I_n;V_h))^2 \times \mathbb P_{k-3}(I_n;Q_h)$.
\end{problem}

\begin{remark}
\begin{itemize}
\item In Problem \ref{Prob:DLSNGC} the variational equation \eqref{Eq:DLSNGC:3} is 
combined with the collocation condition \eqref{Eq:DLSNGC:2} at the endpoint $t_n$ of 
$I_n$ and the continuity constraints \eqref{DLSNGC:1_1}, \eqref{DLSNGC:1_2}.

 \item By definition \eqref{Eq:DefHIntOp} of the Hermite-type interpolation operator 
$\IH$, we have that $\partial_t^s \IH \vec f(t_n) = \partial_t^s \vec f(t_n)$ and 
$\partial_t \IH \vec g(t_n) = \partial_t \vec g(t_n)$ for $s\in \{0,1\}$ on 
the right-hand side of \eqref{Eq:DLSNGC:2}.

\item The choice of the temporal basis (cf.\ Eqs.~\eqref{eq:CondXi}, that is induced 
by the definition of the Hermite-type quadrature formula \eqref{DLSNGC:1_1} and 
the interpolation operator definition \eqref{Eq:DefHIntOp}, allows a computationally 
cost-efficient implementation of the continuity constraints 
\eqref{DLSNGC:1_1}, \eqref{DLSNGC:1_2}. By these constraints the condition $(\vec 
v_{\tau,h},p_{\tau,h})\in \vec X^k_{\tau,k}\cap ((C^1(\overline I;V_h))^2\times 
C^1(\overline I;Q_h))$ and, thus, the $C^1$ regularity in time of 
$\vec v_{\tau,h}$ and $p_{\tau,h}$ is ensured. 

\item For the initial time interval $I_1$, i.e.\ $n=1$, the continuity constraints 
\eqref{DLSNGC:1_1}, \eqref{DLSNGC:1_2} are a source of trouble since we do not have 
an initial pressure $p(0)$ in the Navier--Stokes system \eqref{eq:navier_stokes}. This 
holds similarly to the case of the cGP($k$) approximation in time; cf.\ 
Remark~\ref{Rem:IniPress}. By the construction of the GCC$^1$($k$) approach and its 
temporal basis (cf.\ Eqs.\ \eqref{eq:time_discrete_ansatz}, \eqref{eq:CondXi}), even 
a spatial approximation of the time derivative of the initial pressure 
$\partial_t p(0)$ is needed now. An initial value for $\partial_t \vec v(0)$ and ist 
approximation can still be computed from the momentum equation 
\eqref{eq:navier_stokes_1}. Remedies for the initial time interval $I_1$ are sketched in 
Remark~\ref{Rem:IniPress}. However, this topic still deserves further research in the 
future. In our numerical convergence study presented in 
Subsec.~\ref{sec:convergence_study} the prescribed solution is used for providing the 
needed initial values. In the numerical study of flow around a cylinder presented in 
Subsec.~\ref{sec:flow_cyl} and \ref{sec:flow_re_100} zero initial values are used. This is done to due to the 
specific problem setting.
\end{itemize}

\end{remark}


\section{Algebraic system of  \texorpdfstring{Galerkin--collocation 
GCC$\boldsymbol{{}^1(3)}$}{Galerkin--collocation GCC¹(3) scheme} discretization in 
time and inf-sup stable finite approximation in space and its Newton linearization}
\label{sec:C1_solution}

In this section we derive the algebraic formulation of Problem~\ref{Prob:DLSNGC}. 
The Newton method is applied for solving the resulting nonlinear system of equations. The 
Newton linearization is also developed here. To simplify the notation we restrict 
ourselves to the polynomial degree $k=3$ for the discrete spaces $ (\mathbb 
P_k(I_n;V_h))^2 
\times \mathbb P_k(I_n;Q_h)$. The choice $k=3$ is also used for the numerical experiments 
presented in Sec.~\ref{Sec:NumExp}. To derive the algebraic form of 
Problem~\ref{Prob:DLSNGC}, a Rothe type approach is applied to the system 
\eqref{eq:navier_stokes} by studying firstly in Subsec.~\ref{Subsec:Semdisc} the GCC$^1$(3) discretization in time of the system 
\eqref{eq:navier_stokes} along with its Newton linearization and then, doing  
the discretization in space by the Taylor-Hood family in Subsec.~\ref{Subsec:SpaceDisc}. 
The discretization of the Nitsche terms hidden in the forms 
$a_\gamma$ and $b_\gamma$ of Problem~\ref{Prob:DLSNGC} is derived separately in 
Subsec.~\ref{Subsec:Nitsche}

\subsection{Semi-discretization in time by \texorpdfstring{$GCC^1(3)$}{GCC¹(3)} and 
Newton linearization}
\label{Subsec:Semdisc}

Here, the GCC$^1$(3) discretization in time of the Navier--Stokes system 
\eqref{eq:navier_stokes} and its Newton linearization are presented. To simplify the 
presentation and enhance their confirmability, this is only done formally in the Banach 
space and without providing functions spaces. Further, we assume 
homogeneous Dirichlet boundary conditions $\vec{g} = \vec{0}$ on $\Gamma_i$ in 
this subsection. The extension that are necessary for Nitsche's method are 
sketched in Subsec.~\ref{Subsec:Nitsche}

The GCC$^1$(3) discretization in time of \eqref{eq:navier_stokes} reads as follows. 

\begin{problem}[GCC$^1$(3) semidiscretization in time of \eqref{eq:navier_stokes}]
\label{def:GCC}
Let $k\geq 3$. For $n=1, \ldots, N$ and given $(\vec{v}_{\tau}{}_|{}_{I_{n-1}}, 
p_{\tau}{}_|{}_{I_{n-1}}) \in (\mathbb P_k (I_{n-1};{V}_0))^2\times  \mathbb P_k 
(I_{n-1};W)$ find $(\vec{v}_{\tau}{}_|{}_{I_n}, p_{\tau}{}_|{}_{I_n}) \in (\mathbb P_k 
(I_n;{V}_0))^2 \times \mathbb P_k (I_n;Q)$ such that
\begin{align}
\vec v_{\tau}{}_{|I_n} (t_{n-1})  & = \vec v_{\tau}{}_{|I_{n-1}} (t_{n-1})\,,
& \partial_t \vec v_{\tau}{}_{|I_n} (t_{n-1})  & = \partial_t 
\vec v_{\tau,h}{}_{|I_{n-1}} (t_{n-1})\,,
\label{SDLSGC:1_1}
\\[1ex]
p_{\tau}{}_{|I_n} (t_{n-1}) & = p_{\tau}{}_{|I_{n-1}} (t_{n-1})\,,
& \partial_t p_{\tau}{}_{|I_n} (t_{n-1}) & = \partial_t p_{\tau}{}_{|I_{n-1}} 
(t_{n-1})\,,
\label{SDLSGC:1_2}
\end{align}
and
\begin{align}
\label{SDLSGC:1_3}
\partial_t \vec{v}_{\tau}(t_n) + (\vec{v}_{\tau }(t_n) \cdot 
\nabla) \vec{v}_{\tau}(t_n) - \nu \Delta \vec{v}_{\tau, h}(t_n) + 
\nabla p_{\tau}(t_n) & = \force(t_n)\,, \\[1ex]
\label{SDLSGC:1_4}
\nabla \cdot  \vec{v}_{\tau, h}(t_n) & = 0
\end{align}
and, for all $\zeta_{\tau} \in \mathbb P_{0} (I_n;\R)$,
\begin{align}
\label{SDLSGC:1_5}
\int_{I_n} (\partial_t \vec{v}_{\tau} + (\vec{v}_{\tau} \cdot 
\nabla)\vec{v}_{\tau} - \nu \Delta \vec{v}_{\tau} + \nabla p_{\tau})\cdot  
\zeta_{\tau} \d t & = \int_{I_n} \IH \force \cdot \zeta_{\tau} \d t\,, 
\\[1ex]
\label{SDLSGC:1_6}
\int_{I_n} \nabla \cdot  \vec{v}_{\tau} \cdot \zeta_{\tau} \d t & = 0\,.
\end{align}
\end{problem}

The time integrals in \eqref{SDLSGC:1_5} and \eqref{SDLSGC:1_6} can be computed exactly 
by the quadrature rule \eqref{Eq:GLHF} with $k=3$. For the derivation of an algebraic 
formulation, we firstly rewrite Problem \ref{def:GCC} in terms of conditions about the 
coefficient functions of an expansion of the unknown variables 
$(\vec{v}_{\tau}{}_|{}_{I_n}, p_{\tau}{}_|{}_{I_n}) \in (\mathbb P_k  (I_n;{V}))^2 \times 
\mathbb P_k (I_n;W_h)$ in temporal basis functions $\{\hat \xi_l\}_{l=0}^3$ of $\mathbb 
P_3 (\hat I;\R)$. With the notation $\vec v_{\tau}=(v_{\tau,1},v_{\tau,2})$, such an 
expansion reads as  
\begin{align}
\label{eq:time_discrete_ansatz}
{v}_{\tau,i|I_n}(\mat x, t)
&=
\sum_{l=0}^{3} {v}_{n,l,i}
(\mat{x})\xi_{l}(t)\,,\;\; \text{for } i\in \{1,2\}\,,
&
p_{\tau|I_n}(\mat x, t)
&=
\sum_{l=0}^{3}
p_{n,l}
(\mat{x})\xi_{l}(t)\,,
\end{align}
with coefficient functions $\vec{v}_{n,l} = (v_{n,l,1},v_{n,l,2})\in \vec V_0$ and 
$p_{n,l} \in Q$ and $t\in I_n$. We define the Hermite-type basis $\{\hat \xi_l\}_{l=0}^3$ 
of $\mathbb P_3 (\hat I;\R)$ on the reference time 
interval $\hat I = [0,1]$ by 
\begin{equation}
\label{eq:CondXi}
\begin{aligned}
\hat{\xi_{0}}(0) &= 1\,,
& \hspace{1em}
\hat{\xi_{0}}(1) &= 0\,,
& \hspace{1em}
\partial_t\hat{\xi_{0}}(0) &= 0\,,
& \hspace{1em}
\partial_t\hat{\xi_{0}}(1) &= 0\,,
\\
\hat{\xi_{1}}(0) &= 0\,,
&
\hat{\xi_{1}}(1) &= 0\,,
&
\partial_t\hat{\xi_{1}}(0) &= 1\,,
&
\partial_t\hat{\xi_{1}}(1) &= 0\,,
\\
\hat{\xi_{2}}(0) &= 0\,,
&
\hat{\xi_{2}}(1) &= 1\,,
&
\partial_t\hat{\xi_{2}}(0) &= 0\,,
&
\partial_t\hat{\xi_{2}}(1) &= 0\,,
\\
\hat{\xi_{3}}(0) &= 0\,,
&
\hat{\xi_{3}}(1) &= 0\,,
&
\partial_t\hat{\xi_{3}}(0) &= 0\,,
&
\partial_t\hat{\xi_{3}}(1) &= 1\,.
\end{aligned}
\end{equation}
These conditions yield basis functions of $\mathbb P_3(\hat I;\R)$ that are given by  
\begin{align*}
\hat{\xi_{0}} (\hat t)&= 1 - 3 \hat t^{\,2} + 2 \hat t^{\,3}\,, &
\hat{\xi_{1}} (\hat t) &v = \hat t - 2 \hat t^{\,2} + \hat t^{\,3}\,, &
\hat{\xi_{2}} (\hat t)&= 3 \hat t^{\,2} - 2 \hat t^{\,3}\,, &
\hat{\xi_{3}} (\hat t) &= -\hat t^{\,2} + \hat t^{\,3}
\end{align*}
for $\hat t \in [0,1]$. We note that by \eqref{eq:CondXi} the expansions in 
\eqref{eq:time_discrete_ansatz} thus comprise the function values and time derivatives of 
$\vec{v}_{\tau}{}_|{}_{I_n}$ and $p_{\tau}{}_|{}_{I_n}$ at $t_{n-1}$ and $t_n$.  By 
\eqref{eq:CondXi}, the Hermite-type interpolation operator $\IH$ defined by 
\eqref{Def:I_n} and \eqref{Eq:DefHIntOp} then admits the explicit representation 
\begin{equation}
\label{eq:hermite_interpolation}
\begin{alignedat}{3}
\vec{g}_{\tau}:=I_{\tau|I_n}\vec{g}(t) = 
\sum_{s=0}^{1} \tau_n^s \hat \xi_s  (0)
\underbrace{\partial_t^s \vec{g}_|{}_{I_n}(t_{n-1})}_{=:\vec{g}_{n,s}}
&&+ \sum_{s=0}^{1} \tau_n^s \hat \xi_{s+2} (1)
\underbrace{\partial_t^s \vec{g}|_{I_n}(t_{n})}_{=:\vec{g}_{n, s+2}}\,.
\end{alignedat}
\end{equation}
In terms of the expansions \eqref{eq:time_discrete_ansatz} along with \eqref{eq:CondXi} we recast the conditions \eqref{SDLSGC:1_1} and \eqref{SDLSGC:1_2}  as
\begin{equation}
\label{eq:time_discrete_collocation_1}
\begin{aligned}
\vec{v}_{n,0} &= \vec{v}_{n-1,2}\,,
&
\vec{v}_{n,1} &= \vec{v}_{n-1,3}\,,
&
p_{n,0} &= p_{n-1,2}\,,
&
p_{n,1} &= p_{n-1,3}\,.
\end{aligned}
\end{equation}
Therefore, on each time interval $I_n$ we obtain four unknown coefficient functions from 
the previous time interval $I_{n-1}$ which amounts to computationally cheap copies of 
vectors on the fully discrete level. Next, integrating analytically the variational 
conditions \eqref{SDLSGC:1_5} and \eqref{SDLSGC:1_6} with the representations 
\eqref{eq:time_discrete_ansatz} of the unknowns implies that
\begin{equation}
\label{eq:time_discrete_variational_3a}
\begin{aligned}
\vec{v}_{n,0} & + \vec{v}_{n,2}
+
\tau_n \biggl(
\frac{13}{35} \nabla \vec{v}_{n,0} \vec{v}_{n,0}
+ \frac{11}{210} \nabla \vec{v}_{n,0} \vec{v}_{n,1}
+ \frac{9}{70} \nabla \vec{v}_{n,0} \vec{v}_{n,2}
- \frac{13}{420} \nabla \vec{v}_{n,0} \vec{v}_{n,3}
\\
& + \frac{11}{210} \nabla \vec{v}_{n,1} \vec{v}_{n,0}
+ \frac{1}{105} \nabla \vec{v}_{n,1} \vec{v}_{n,1}
+ \frac{13}{420} \nabla \vec{v}_{n,1} \vec{v}_{n,2}
- \frac{1}{140} \nabla \vec{v}_{n,1} \vec{v}_{n,3}
+ \frac{9}{70} \nabla \vec{v}_{n,2} \vec{v}_{n,0}
\\
& + \frac{13}{420} \nabla \vec{v}_{n,2} \vec{v}_{n,1}
+ \frac{13}{35} \nabla \vec{v}_{n,2} \vec{v}_{n,2}
- \frac{11}{210} \nabla \vec{v}_{n,2} \vec{v}_{n,3}
- \frac{13}{420} \nabla \vec{v}_{n,3} \vec{v}_{n,0}
- \frac{1}{140} \nabla \vec{v}_{n,3} \vec{v}_{n,1}
\\
& - \frac{11}{210} \nabla \vec{v}_{n,3} \vec{v}_{n,2}
+ \frac{1}{105} \nabla \vec{v}_{n,3} \vec{v}_{n,0}
\biggr)
- \tau_n \nu \biggl(
\frac{1}{2} \Delta \vec{v}_{n,0}
+ \frac{1}{12} \Delta \vec{v}_{n,1}
+ \frac{1}{2} \Delta \vec{v}_{n,2}
- \frac{1}{12} \Delta \vec{v}_{n,3}
\biggr)
\\
&+ \tau_n \biggl(
\frac{1}{2} \nabla p_{n, 0}
+ \frac{1}{12} \nabla p_{n, 1}
+ \frac{1}{2} \nabla p_{n, 2}
- \frac{1}{12} \nabla p_{n, 3}
\biggr)
=
\tau_n \biggl(
\frac{1}{2} \nabla f_{n, 0}
+ \frac{1}{12} \nabla f_{n, 1}
+ \frac{1}{2} \nabla f_{n, 2}
- \frac{1}{12} \nabla f_{n, 3}
\biggr)
\end{aligned}
\end{equation}
and
\begin{equation*}
\tau_n \biggl(
\frac{1}{2} \nabla \cdot \vec{v}_{n,0}
+ \frac{1}{12} \nabla \cdot \vec{v}_{n,1}
+ \frac{1}{2} \nabla \cdot \vec{v}_{n,2}
- \frac{1}{12} \nabla \cdot \vec{v}_{n,3}
\biggr)
= 0 \,.
\end{equation*}

We note that in order to keep the notation as short as possible, we omit the brackets in 
$\left(\nabla \vec{v}_{n,i}\right)\vec{v}_{n,j}$ and assume that the differential 
operator just acts on the vector next to it only such that $\left(\nabla 
\vec{v}_{n,i}\right)\vec{v}_{n,j} = \nabla \vec{v}_{n,i}\vec{v}_{n,j}$. Finally, by \eqref{eq:time_discrete_ansatz} along with \eqref{eq:CondXi} the 
collocation conditions\eqref{SDLSGC:1_3} and \eqref{SDLSGC:1_4} read as 
\begin{align}
\frac{1}{\tau_n} \vec{v}_{n,3}
\,\,+
\nabla \vec{v}_{n,2} \vec{v}_{n,2}
\,\,-
\nu \Delta \vec{v}_{n,2}
+
\nabla p_{n,2}
&=
\force_{n,2}\,,
\label{eq:time_discrete_collocation_2a}
\\
\nabla \cdot \vec{v}_{n, 2}
&=
0\,,
\label{eq:time_discrete_collocation_2b}
\end{align}

On each subinterval $I_n$,  Eqs.~\eqref{eq:time_discrete_collocation_1} to \eqref{eq:time_discrete_collocation_2b} form a nonlinear system in the Banach space. To solve this system of nonlinear equations (after an additional discretization in space; cf. Subsec.~\ref{Subsec:SpaceDisc}) we use an inexact Newton method. To enhance the region of convergence, a line search approach damping the length of the Newton step was implemented. Moreover, the dogleg method (cf., e.g., \cite{P08}), which belongs to the class of trust-region methods, was tested. The advantage of the latter one is that also the search direction, and not just its length, is adapted. For both schemes, the Jacobian matrix related to Eqs.\ \eqref{eq:time_discrete_collocation_2a}, \eqref{eq:time_discrete_collocation_2b} has to be computed. In the dogleg method, multiple products of the Jacobian matrix with vectors are needed further. Since the Jacobian is stored as a sparse matrix, this matrix--vector products can be computed at low computational costs.
From the theoretical point of view both methods yield superlinear convergence. For our numerical examples of  Sec.~\ref{Sec:NumExp} both modifications of the standard Newton method lead to comparable results. We did not observe any convergence problems. Since the Galerkin--collocation method is of higher order in time and, thereby, aims at using large time step sizes for reasonable accuracy, fixed-point iterations with first-order convergence only and modifications of them, like the L-scheme (cf.~\cite{LR16}), were not considered here. If convergence problems arise in future applications, a combination of fixed-point methods and Newton iteration are an option to ensure convergence.

For the sake of completeness, the (non-damped) Newton iteration for the system \eqref{eq:time_discrete_collocation_1} to \eqref{eq:time_discrete_collocation_2b} in the Banach space is briefly sketched here. For this, we let 
\begin{equation*}
\vec x \coloneqq (\vec{v}_{n,2}, p_{n,2}, \vec{v}_{n,3}, p_{n, 3})^\top
\end{equation*}
denote the vector of remaining unknown coefficient functions of the expansions \eqref{eq:time_discrete_ansatz} under the identities \eqref{eq:time_discrete_collocation_1} . Further, we subtract the right-hand sides of Eqs.\ \eqref{eq:time_discrete_variational_3a} to  \eqref{eq:time_discrete_collocation_2b} from its left-hand sides respectively and denote the resulting left-hand side functions by $\vec{q}(\vec{x})$ with its components
$\{\vec{q}_1(\vec{x})$, $q_2(\vec{x})$, $\vec{q}_3(\vec{x})$, $q_4(\vec{x})\}$. Then, the Newton iteration reads as 
\begin{equation}
\label{eq:Newton_2}
\mat{J}_{\vec{q},{\vec{x}^{k}}} \delta \vec{x}^{k+1}
= - \vec{q}(\vec{x}^k) \,.
\end{equation}
for the correction $\delta \vec{x}^{k+1}: = \vec{x}^{k+1} - \vec{x}^k$ with the 
directional derivative along $\delta \vec{x}^{k+1}$ at $\vec{x}^k$ given by 
\begin{equation*}
	\mat{J}_{\vec{q},{\vec{x}^{k}}} \delta \vec{x}^{k+1}
	=
	\lim\limits_{\epsilon \rightarrow 0}
	\frac{1}{\epsilon}
	\Bigl(
	\vec{q}(\vec{x}^k + \epsilon \delta \vec{x}^{k})
	-
	\vec{q}(\vec{x}^k)
	\Bigr)\,.
\end{equation*}
Introducing the abbreviations 
\begin{align*}
	\vec{a} (\vec{x}) &= 
	\left(
	\frac{9}{70} \vec{v}_{n,0}
	+\frac{13}{420} \vec{v}_{n,1}
	+\frac{13}{35} \vec{v}_{n,2}
	-\frac{11}{210} \vec{v}_{n,3}
	\right)\,,
	\\
	\vec{b} (\vec{x}) &=
	\left(\frac{-13}{420} \vec{v}_{n,0} 
	-\frac{1}{140} \vec{v}_{n,1}
	-\frac{11}{210} \vec{v}_{n,2}
	+\frac{1}{105} \vec{v}_{n,3}
	\right) 
\end{align*}
and defining
\begin{equation*}
	\begin{aligned}
	& \begin{aligned}
	\vec{f}_1(\delta \vec{v}_{n,2}^{k+1}) \coloneqq
	\delta \vec{v}_{n,2} ^{k+1}
	-\frac{\tau_n \nu}{2} \Delta \delta \vec{v}_{n,2}^{k+1}
	+ \tau_n \nabla \delta \vec{v}_{n,2}^{k+1} \vec{a}(\vec{x}^k)
	+ \tau_n \nabla \vec{a}(\vec{x}^k) \delta \vec{v}_{n,2}^{k+1} \,,
	\end{aligned}
	\\
	& \begin{aligned}
	\vec{f}_2(\delta \vec{v}_{n,3}^{k+1}) \coloneqq
	 \frac{\tau_n \nu}{12} \Delta \delta \vec{v}_{n,3}^{k+1}
	+ \tau_n \nabla \delta \vec{v}_{n,3}^{k+1} \vec{b}(\vec{x}^k)
	+ \tau_n \nabla \vec{b}(\vec{x}^k) \delta \vec{v}_{n,3}^{k+1} \,,
	\end{aligned}
	\\
	& \vec{f}_3(\delta \vec{v}_{n,2}^{k+1}) \coloneqq
	\nu \Delta \delta\vec{v}_{n,2}^{k+1}
	+ \nabla \vec{v}_{n,2}^k \delta\vec{v}_{n,2}^{k+1}
	+ \nabla \delta\vec{v}_{n,2}^{k+1} \vec{v}_{n,2}^k \,,
	\quad
	\vec{f}_4(\delta \vec{v}_{n,3}^{k+1}) \coloneqq
	\frac{1}{\tau_n} \delta\vec{v}_{n,3}^{k+1} \,,
	\\
	& \begin{aligned}
	\vec{b}_1 (\delta \vec{v}_{n,2}^{k+1}) &\coloneqq
	- \frac{\tau_n}{2} \nabla \cdot \delta \vec{v}_{n,2}^{k+1} \,,
	&
	\vec{b}_2 (\delta \vec{v}_{n,3}^{k+1}) &\coloneqq
	- \frac{\tau_n}{12} \nabla \cdot \delta \vec{v}_{n,3}^{k+1} \,,
	&
	\vec{b}_3 (\delta \vec{v}_{n,3}^{k+1}) &\coloneqq
	- \nabla \cdot \delta \vec{v}_{n,3}^{k+1} \,,
	\\
	\vec{b}_1^\top (\delta p_{n,2}^{k+1}) & \coloneqq
	- \frac{\tau_n}{2} \nabla \delta p_{n,2}^{k+1} \,,
	&
	\vec{b}_2^\top (\delta p_{n,3}^{k+1}) &\coloneqq
	- \frac{\tau_n}{12} \nabla \delta p_{n,3}^{k+1} \,,
	&
	\vec{b}_3^\top (\delta p_{n,3}^{k+1}) &\coloneqq
	- \nabla \delta p_{n,3}^{k+1} \,,
	\end{aligned}
	\end{aligned}
\end{equation*}
we recast the system \eqref{eq:Newton_2} as 
\begin{subequations}
	\label{Eq:NewSys}
\begin{align}
	\vec{f}_1 (\delta \vec{v}_{n,2}^{k+1}) + \vec{f}_2 (\delta \vec{v}_{n,3}^{k+1}) 
	+ \vec{b}_1^\top (\delta p_{n,2}^{k+1}) + \vec{b}_2^\top (\delta p_{n,3}^{k+1}) &= 
- \vec{q}_1 (\vec{x}^k)\,,
	\label{eq:Newton_3}
	\\[1ex]
	- \vec{b}_1 (\delta \vec{v}_{n,2}^{k+1}) + \vec{b}_2 (\delta \vec{v}_{n,3}^{k+1}) 
&= - q_2 (\vec{x}^k)\,,
	\label{eq:Newton_4}
	\\[1ex]
	\vec{f}_3 (\delta \vec{v}_{n,2}^{k+1}) + \vec{f}_4 (\delta \vec{v}_{n,3}^{k+1})
	 + \vec{b}_3^\top (\delta p_{n,3}^{k+1}) &= - \vec{q}_3 (\vec{x}^k)\,,
	\label{eq:Newton_5}
	\\[1ex]
	-\vec{b}_3 (\delta \vec{v}_{n,3}^{k+1}) &= - q_4 (\vec{x}^k)\,.
	\label{eq:Newton_6}
\end{align}
\end{subequations}
In weak form, system \eqref{Eq:NewSys} leads to the following problem to be solved 
in each Newton step. 

\begin{problem}[Newton iteration of GCC$^1$(3) time discretization]
\label{Prob:ContNewton}
Find corrections $(\delta \vec{v}_{n,2}^{k+1}, \delta \vec{v}_{n,3}^{k+1}) \in \mat{V}_{{0}}^2$ and $(\delta p_{n,2}^{k+1}, \delta p_{n,3}^{k+1})\in Q^2$ such that for all $\vec{\psi} \in \mat{V}_{{0}}^2$ and $\xi \in Q^2$ there holds that 
\begin{subequations}
\label{eq:Newton_continuous_0}
\begin{align}
&\begin{aligned}
\mathcal{F}_1 (\delta \vec{v}_{n,2}^{k+1}, \vec{\psi}) 
&+ \mathcal{F}_2 (\delta \vec{v}_{n,3}^{k+1}, \vec{\psi})
 + \mathcal{B}_1^\top (\delta p_{n,2}^{k+1}, \vec{\psi})
+ \mathcal{B}_2^\top (\delta p_{n,3}^{k+1}, \vec{\psi})
\\[1ex]
&= - \Bigl<
\vec{q}_1 (\vec{x}^k)
\,, \vec{\psi} \Bigr>
+
\frac{\tau_n \nu}{2} \Bigl<
\partial_{\vec{n}}\delta\vec{v}_{n, 2} (\vec{x}^k)
\,, \vec{\psi} \Bigr>_{\Gamma_o}
-
\frac{\tau_n \nu}{12} \Bigl<
\partial_{\vec{n}}\delta\vec{v}_{n, 3} (\vec{x}^k)
\,, \vec{\psi} \Bigr>_{\Gamma_o} \,,
\label{eq:Newton_continuous_1}
\end{aligned}
\\[1ex]
&- \mathcal{B}_1 (\delta \vec{v}_{n,2}^{k+1}, \xi) 
- \mathcal{B}_2 (\delta \vec{v}_{n,3}^{k+1}, \xi) 
=
- \Bigl<
q_2 (\vec{x}^k)
\,, \xi \Bigr>\,,
\label{eq:Newton_continuous_2}
\\[1ex]
&\mathcal{F}_3 (\delta \vec{v}_{n,2}^{k+1},  \vec{\psi}) 
+ \mathcal{F}_4 (\delta \vec{v}_{n,3}^{k+1},  \vec{\psi})
+ \mathcal{B}_3^\top (\delta p_{n,3}^{k+1},  \vec{\psi})
 = 
- \Bigl<
\vec{q}_3 (\vec{x}^k)
\,, \vec{\psi} \Bigr>
- \nu \Bigl<
\partial_{\vec{n}}\delta\vec{v}_{n, 2} (\vec{x}^k)
\,, \vec{\psi} \Bigr>_{\Gamma_o} \,,
\label{eq:Newton_continuous_3}
\\[1ex]
&-\mathcal{B}_3 (\delta \vec{v}_{n,3}^{k+1}, \xi) 
 =
- \Bigl<
q_4 (\vec{x}^k)
\,, \xi \Bigr>\,,
\label{eq:Newton_continuous_4}
\end{align}
\end{subequations}
with $\partial_{\vec{n}} \vec{w} = \nabla \vec{w} \cdot \vec{n}$ and
\begin{equation*}
\begin{aligned}
	& \begin{multlined}[c][\textwidth]
	\mathcal{F}_1(\delta \vec{v}_{n,2}^{k+1}, \vec{\psi}) \coloneqq
	\Bigl< \delta \vec{v}_{n,2} ^{k+1} , \vec{\psi} \Bigr>
	+ \frac{\tau_n \nu}{2} \Bigl<
	\nabla \delta \vec{v}_{n,2}^{k+1}
	, \nabla \vec{\psi} \Bigr>
	+ \tau_n \Bigl<
	\nabla \delta \vec{v}_{n,2}^{k+1} \vec{a}(\vec{x}^k)
	, \vec{\psi} \Bigr>
	+ \tau_n \Bigl<
	\nabla \vec{a}(\vec{x}^k) \delta \vec{v}_{n,2}^{k+1}
	, \vec{\psi} \Bigr>\,,
	\end{multlined}
	\\[0.5ex]
	& \begin{multlined}[c][\textwidth]
	 \mathcal{F}_2(\delta \vec{v}_{n,3}^{k+1}, \vec{\psi}) \coloneqq
	- \frac{\tau_n \nu}{12} \Bigl<
	\nabla \delta \vec{v}_{n,3}^{k+1}
	, \nabla \vec{\psi} \Bigr>
	+ \tau_n \Bigl<
	\nabla \delta \vec{v}_{n,3}^{k+1} \vec{b}(\vec{x}^k)
	, \vec{\psi} \Bigr>
	+ \tau_n \Bigl<
	\nabla \vec{b}(\vec{x}^k) \delta \vec{v}_{n,3}^{k+1}
	, \vec{\psi} \Bigr>\,,
	\end{multlined}
	\\[0.5ex]
	& \mathcal{F}_3(\delta \vec{v}_{n,2}^{k+1}, \vec{\psi}) \coloneqq
	\nu \Bigl<
	\nabla \delta\vec{v}_{n,2}^{k+1}
	, \nabla \vec{\psi} \Bigr>
	+ \Bigl<
	\nabla \vec{v}_{n,2}^k \delta\vec{v}_{n,2}^{k+1}
	, \vec{\psi} \Bigr>
	+ \Bigl<
	\nabla \delta\vec{v}_{n,2}^{k+1} \vec{v}_{n,2}^k
	, \vec{\psi} \Bigr> \,,
	\\[0.5ex]
	&\mathcal{F}_4(\delta \vec{v}_{n,3}^{k+1}, \vec{\psi}) \coloneqq
	\frac{1}{\tau_n} \Bigl<
	\delta \vec{v}_{n,3} ^{k+1} , \vec{\psi}
	\Bigr> \,,
	\\[0.5ex]
\end{aligned}
\end{equation*}
and 
\begin{equation}
\label{eq:bilinear_forms}
\begin{aligned}
	& \mathcal{B}_1 (\delta \vec{v}_{n,2}^{k+1}, \xi) \coloneqq
	- \frac{\tau_n}{2} \Bigl<
	\nabla \cdot \delta \vec{v}_{n,2}^{k+1}
	, \xi \Bigr>
	\,, \quad
	\mathcal{B}_1^\top (\delta p_{n,2}^{k+1}, \vec{\psi}) \coloneqq
	\frac{\tau_n}{2} \Bigl<
	\delta p_{n,2}^{k+1}
	, \nabla \cdot \vec{\psi} \Bigr>\,,
	\\[0.5ex]
	& \mathcal{B}_2 (\delta \vec{v}_{n,2}^{k+1}, \xi) \coloneqq
	\frac{\tau_n}{12} \Bigl<
	\nabla \cdot \delta \vec{v}_{n,3}^{k+1}
	, \xi \Bigr>
	\,, \quad
	\mathcal{B}_2^\top (\delta p_{n,3}^{k+1}, \vec{\psi}) \coloneqq
	- \frac{\tau_n}{12} \Bigl<
	\delta p_{n,3}^{k+1}
	, \nabla \cdot \vec{\psi} \Bigr>\,,
	\\[0.5ex]
	& \mathcal{B}_3 (\delta \vec{v}_{n,2}^{k+1}, \xi) \coloneqq
	- \Bigl<
	\nabla \cdot \delta \vec{v}_{n,3}^{k+1}
	, \xi \Bigr>
	\,, \quad
	\mathcal{B}_3^\top (\delta p_{n,3}^{k+1}, \vec{\psi}) \coloneqq
	\Bigl<
	\delta p_{n,3}^{k+1}
	, \nabla \cdot \vec{\psi} \Bigr>\,.
\end{aligned}
\end{equation}
\end{problem}


\subsection{Fully discrete system with inf-sup stable elements}
\label{Subsec:SpaceDisc}

In this subsection we briefly present the discretization in space of the system 
\eqref{eq:Newton_continuous_0} of Problem~\ref{Prob:ContNewton} in the pair $V_h^0\times 
Q_h$ of inf-sup stable finite element spaces. For our computations presented in 
Sec.~\ref{Sec:NumExp} we used the $\mathcal Q_p$--$\mathcal Q_{p-1}, p \geq 2$ pair of the 
well-known Taylor--Hood family. Due to their inf-sup stability a stabilization of the 
discretization is not required, as long as the Reynolds number of the fluid flow is 
assumed to be small such that no convection-dominance occurs. In the case of higher 
Reynolds numbers, an additional stabilization of the discretization becomes 
indispensible; cf.~\cite[Sec.~5.3 and 5.4]{J16} and the references therein. However, this 
is beyond the scope of interest in this work and left as a work for the future.

For the space discretization, let $\{\vec{\psi}_j\}_{j=1}^J \subset 
\vec{V}_{h}^0$ and $\{\phi_m\}_{m=1}^M \subset Q_h$ denote a nodal 
Lagrangian basis of $\vec{V}_h^0$ and $Q_h$, respectively. Then, the fully discrete 
unknows admit the representations 
\begin{equation*}
\begin{aligned}
\vec{v}_{\tau,h}{}_{|I_n}(\mat x, t)
&=
\sum_{l=0}^{3} \sum_{j=1}^{J}
v_{n,l,j} \vec{\psi}_j
(\mat{x})\xi_{l}(t)\,,
&
p_{\tau,h}{}_{|I_n}(\mat x, t)
&=
\sum_{l=0}^{3} \sum_{m=1}^{M}
p_{n,l,m} \phi_m
(\mat{x})\xi_{l}(t)
\end{aligned}
\end{equation*}
for $(\mat x,t)\in \Omega \times \overline{I_n}$ with the unknown coefficient vector $\vec v_{n,l}:=(v_{n,l,j})_{j=1}^J\in \R^J$ and $\vec p_{n,l}:=(p_{n,l,m})_{m=1}^{M}\in \R^M$, for $l=0,\ldots, 3$, as the degrees of freedom. Next, we define
\begin{subequations}
\label{Def:F_iB_i}
\begin{align}
	\mat{F}_1 & \coloneqq \left( \mathcal{F}_1(\vec{\psi}_i, \vec{\psi}_j) \right)_{i,j=1}^J \,,
	& 
	\mat{F}_2 & \coloneqq \left( \mathcal{F}_2(\vec{\psi}_i, \vec{\psi}_j) \right)_{i,j=1}^J \,,
	&
	\mat{F}_3 & \coloneqq \left( \mathcal{F}_3(\vec{\psi}_i, \vec{\psi}_j) \right)_{i,j=1}^J \,,
	\\
	\mat{F}_4 & \coloneqq \left( \mathcal{F}_4(\vec{\psi}_i, \vec{\psi}_j) \right)_{i,j=1}^J \,,
	&
	\mat{B}_1 & \coloneqq \left( \mathcal{B}_1(\vec{\psi}_i, \phi_l) \right)_{i, m  = 
1}^{J,M} \,,
	&
	\mat{B}_1^\top & \coloneqq \left( \mathcal{B}_1^\top(\vec{\psi}_i, \vec{\psi}_j) \right)_{i, j  = 1}^J \,,
	\\
	\mat{B}_2 & \coloneqq \left( \mathcal{B}_2(\vec{\psi}_i, \phi_m) \right)_{i, m  = 1}^{J,M} \,,
	&
	\mat{B}_2^\top & \coloneqq \left( \mathcal{B}_2^\top(\vec{\psi}_i, \vec{\phi}_j) \right)_{i, j  = 1}^J \,,
	&
	\mat{B}_3 & \coloneqq \left( \mathcal{B}_3(\vec{\psi}_i, \phi_m) \right)_{i, m  = 1}^{J,M} \,,
	\\
	\mat{B}_3^\top & \coloneqq \left( \mathcal{B}_3^\top(\vec{\psi}_i, \vec{\psi}_j) \right)_{i, j  = 1}^J \,.
\end{align}
\end{subequations}

Solving Problem \ref{Prob:ContNewton} in the finite dimensional subspaces $\vec V_h^0$ 
and $Q_h$ of $\vec V_0$ and $Q$, respectively, leads to the following  problem to be 
solved within each Newton iteration of a time step. 

\begin{problem}[Newton iteration of GCC$^1$(3) time discretization and inf-sup stable elements for space discretization]
\label{Prob:NewtonDiscret}	
Find corrections $(\delta \vec{v}_{n,2}^{k+1}, \delta \vec{v}_{n,3}^{k+1}) \in \R^{2J}$ and $(\delta \vec{p}_{n,2}^{k+1}, \delta \vec{p}_{n,2}^{k+1})\in \R^{2M}$ such that
\begin{equation}
\label{eq:resulting_problem}
	\mat{S} \delta \vec{x}^{k+1} = \vec{d}(\vec{x}^k) \,,
\end{equation}
where $\delta\vec{x}^{k+1} \coloneqq (\delta\vec{v}_{n,2}^{k+1}, 
\delta\vec{p}_{n,2}^{k+1}, \delta\vec{v}_{n,3}^{k+1}, \delta\vec{p}_{n, 3}^{k+1})^\top 
\in \R^{2(J+M)}$ denotes the vector of Newton corrections for the degrees of freedom and 
$\vec{d}(\vec{x}^k)$ is the fully discrete counterpart of the terms on the right-hand side 
of \eqref{eq:Newton_continuous_0}. The block system matrix $\mat S$ in 
\eqref{eq:resulting_problem} is given by
\begin{equation}
\label{eq:newton_system_matrix}
	\mat{S} =
	\begin{pmatrix}
	\mat{F}_1 & \mat{B}_1^\top & \mat{F}_2 & \mat{B}_2^\top \\[1ex]
	-\mat{B}_1 & \mat{0} & -\mat{B}_2 & \mat{0} \\[1ex]
	\mat{F}_3 & \mat{B}_3^\top & \mat{F}_4 & \mat{0} \\[1ex]
	-\mat{B}_3 & \mat{0} & \mat{0} & \mat{0}
	\end{pmatrix} .
\end{equation}
\end{problem}

\begin{figure}[!htb]
	\centering
	\includegraphics[width=0.6\textwidth,keepaspectratio]
	{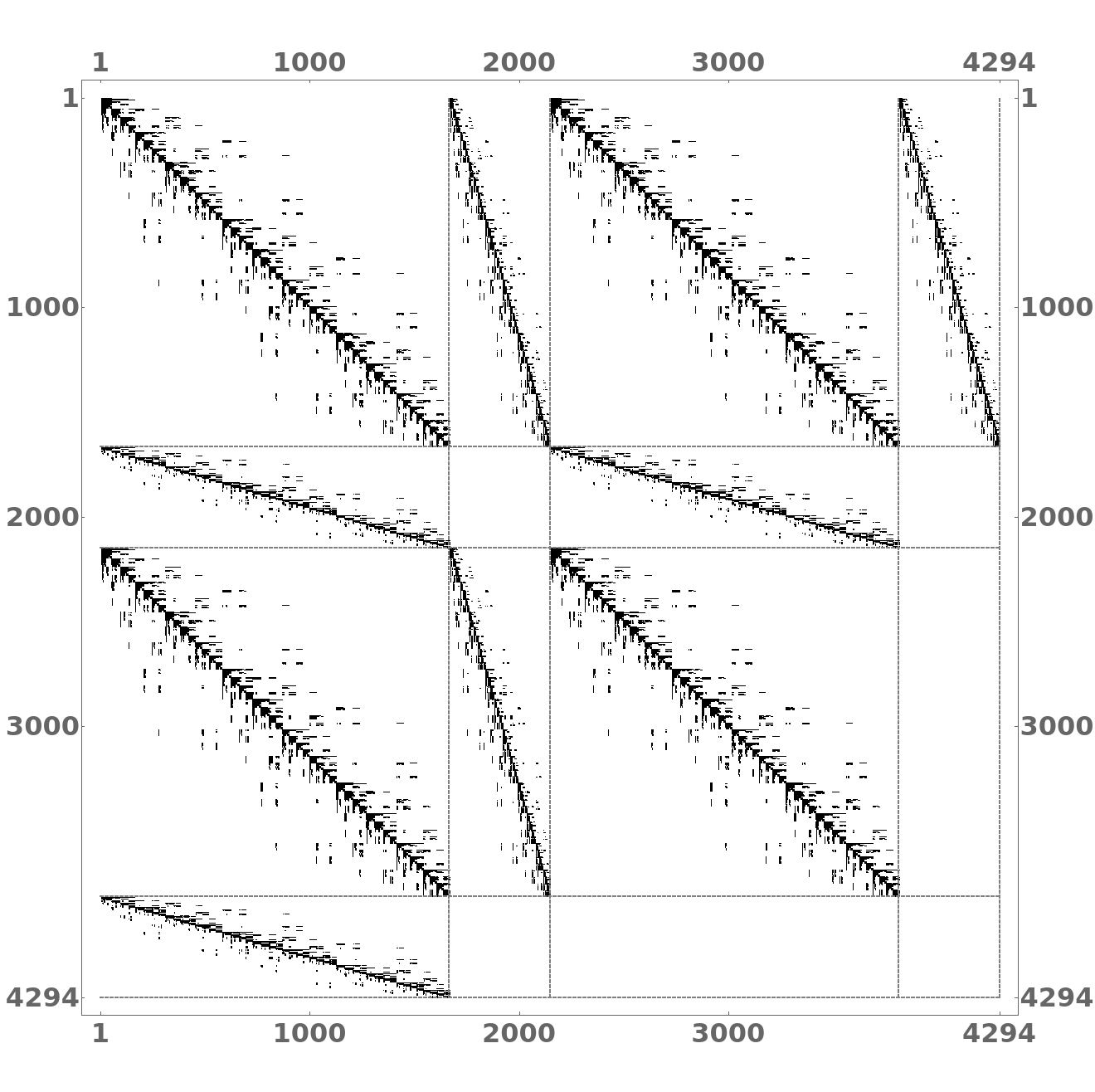}
	\caption{Sparsity pattern of $\vec S$ in \eqref{eq:newton_system_matrix} for the 
discrete problem of Subsec.~\ref{sec:convergence_study}, using $\mathcal 
Q_4$--$\mathcal Q_3$ Taylor-Hood elements.}
	\label{fig:sparsity_pattern}
\end{figure}

Since we use the family of inf-sup stable Taylor-Hood element here, the resulting 
system matrix \eqref{eq:newton_system_matrix} comprises non-quadratic sub-matrices 
$\mat{B}$. The sparsity pattern of $\vec S$ is illustrated in 
Fig.~\ref{fig:sparsity_pattern}. The system matrix $\mat{S}$ consists of three submatrices 
$\mat{S}_i$ of the common structure
\begin{equation}
	\label{eq:newton_system_matrix_1}
	\mat{S}_i =
	\begin{pmatrix}
	\mat{F}_i & \mat{B}_i^\top \\[1ex]
	-\mat{B}_i & \mat{0}
	\end{pmatrix} \,,
\end{equation}
and an additional block of $\mat{F}_4$ together with blocks of zero entries. Due to 
the collocation conditions at the final time points $t_n$ of the subinterval $I_n$ the 
matrices  $\mat{B}_4^\top$ and a $-\mat{B}_4$ do not arise in the right lower block of 
$\vec S$ such that in \eqref{eq:newton_system_matrix} a sparser matrix structure is 
obtained compared to a pure variational approach. In order to solve 
\cref{eq:resulting_problem} we use a (parallel) GMRES solver with a (preliminary) 
block preconditioner that is motivated by an approach presented in \cite{ESW14}. In 
\eqref{eq:newton_system_matrix}, we consider each of the three submatrices $S_i$ as an 
uncoupled block of the structure \eqref{eq:newton_system_matrix_1}. For each of 
these blocks we then use a Schur complement preconditioner with an approximation of 
the mass matrix of the pressure variable. This results in reasonable numbers of 
iterations for the GMRES solver for two-dimensional problems of medium 
size but is far from being acceptable for three-dimensional or large scale problems. The 
design of a more efficient and robust preconditioner that is tailored to the specific 
structure of the matrix $\vec S$ in \eqref{eq:newton_system_matrix} or a multigrid 
approach remains as a work for the future. 

\subsection{Nitsche's method for boundary conditions of Dirichlet type}
\label{Subsec:Nitsche}

In this subsection we briefly present the modifications to be made 
in Problem~\ref{Prob:ContNewton} and \ref{Prob:NewtonDiscret}, respectively, for the 
application of Nitsche's method to enforce Dirichlet boundary conditions; cf.\ 
Problem~\ref{Prob:DLSNGC}. In contrast to Problem~\ref{Prob:ContNewton} and 
\ref{Prob:NewtonDiscret}, the  Dirichlet boundary conditions are now ensured by augmenting 
the weak formulation with additional line integrals (surface integrals in three space 
dimensions); cf.~Problem~\ref{Prob:DLSNGC}. In the field of computational 
fluid dynamics, Nitsche's method offers appreciable advantages over the 
standard implementation of Dirichlet boundary condition and is particularly well suited 
if complex and dynamic geometries are considered. The geometry can be immersed into an 
underlying computational grid. The Navier-Stokes equations are then solved fulfilling the 
boundary conditions at the intersections between the surface discretization and the grid 
cells; cf.\ Sec.~\ref{Sec:Outlook}.

In contrast to Sec.~\ref{Subsec:Semdisc}, the continuous solution and test space is now 
$\vec X = \vec V \times Q$ instead of $\vec X_0 = \vec V_0 \times Q$. For the weak form of 
the time discrete Newton linearization \eqref{Eq:NewSys} integration by parts then yields 
that 
\begin{subequations}
\label{Eq:IntParts}
\begin{align}
\label{Eq:IntParts_1}
- \nu \bigl< \Delta \delta \vec{v}_{\tau}^{k+1}, \vec{\psi} \bigr>
& = \nu \bigl< \nabla \delta \vec{v}_{\tau}^{k+1}, \nabla \vec{\psi} \bigr>
- \nu \bigl< \partial_{\vec{n}} \delta \vec{v}_{\tau}^{k+1}, \vec{\psi} \bigr>_{\Gamma}\,,\\[1ex]
\label{Eq:IntParts_2}
\bigl< \nabla \delta p_{\tau}^{k+1}, \vec{\psi} \bigr>
& = - \bigl<\delta p_{\tau}^{k+1}, \nabla \cdot \vec{\psi} \bigr>
+ \bigl< \delta p_{\tau}^{k+1} \vec{n}, \vec{\psi} \bigr>_{\Gamma}
\end{align}
\end{subequations}
where $\bigl< \cdot, \cdot \bigr>_{\Gamma}$ denotes the inner product of $L^2(\partial 
\Omega)$ and $\vec g_\tau$ is defined by means of the Hermite-type interpolation 
\eqref{eq:hermite_interpolation}. To preserve the symmetry properties of the continuous 
system, the forms (cf.\ Problem~\ref{Prob:DLSNGC})
\begin{subequations}
\label{eq:Nitsche_0}
\begin{align}
\label{eq:Nitsche_1}
S_v (\vec \psi) & := \nu \bigl< \partial_{\vec{n}} \vec{\psi}, \vec{v}_{\tau}^k + \delta \vec{v}_{\tau}^{k+1} - \vec{g}_{\tau} \bigr>_{\Gamma_D}\,, \\[1ex]
\label{eq:Nitsche_2}
S_p(\xi)  & := \bigl< \xi \vec{n}, \vec{v}_{\tau}^k + \delta \vec{v}_{\tau}^{k+1} - \vec{g}_{\tau} \bigr>_{\Gamma_D} 
\end{align}
\end{subequations}
are added on the right-hand side of \eqref{Eq:IntParts_1} and  
\eqref{Eq:IntParts_2}, respectively, to the viscous and pressure part. Finally, we add 
the penalty terms (cf.\ Problem~\ref{Prob:DLSNGC})
\begin{equation}
\label{eq:Nitsche_3}
P_{\vec g_\tau}(\vec \psi) := 
	\underbrace{
	\frac{\eta_1}{h} \nu \bigl< \vec{v}_{\tau}^k + \delta \vec{v}_{\tau}^{k+1} - \vec{g}_{\tau}, \vec{\psi} \bigr>_{\Gamma_D}
	}_{\text{viscous penalty}} + 
	\underbrace{
	\frac{\eta_2}{h} \bigl< (\vec{v}_{\tau}^k + \delta \vec{v}_{\tau}^{k+1} - \vec{g}_{\tau}) \cdot \vec{n}, \vec{\psi} \cdot \vec{n} \bigr>_{\Gamma_D}
	}_{\text{penalty along $\vec{n}$}} \,,
\end{equation}
For the test space $\vec X =\vec V\times Q$, such that $(\vec \psi, \xi) \in \vec V\times 
Q$, the integrals over the Dirichlet part $\Gamma_D$ of the boundary $\Gamma = 
\Gamma_D\cup \Gamma_o$ no longer vanish. For viscous-dominated flow, the additional terms 
in \eqref{eq:Nitsche_3}  enforce in weak form the boundary condition $\vec{v} - \vec{g} = 
\vec{0}$ on the Dirchet part of the boundary. For convection-dominated flow with a small 
viscosity parameter $\nu$ this enforcement is weakened in the first term on the 
right-hand side of \eqref{eq:Nitsche_3}. However, in the inviscid limit $\nu \rightarrow 
0$, the condition $(\vec{v} - \vec{g}) \cdot \vec{n} = 0$ is still imposed weakly by 
the second of the terms on the right-hand side of \eqref{eq:Nitsche_3} such that the 
normal component of the Dirichlet boundary condition is preserved in the limit case. For  
our computations shown in Sec.~\ref{Sec:NumExp} we put  $\eta_1 = \eta_2 = 35$. For a 
more refined analysis of these parameters we refer to \cite{MSW18,S17}. Changing the sign 
of the symmetric term \eqref{eq:Nitsche_1} generates a non-symmetric formulation. Current 
results (cf.\ \cite{S17,BB16,B12}) show that, on the one hand, this reduces the 
sensitivity with respect of the choice of the penalty parameters and might even allow for 
a parameter free penalty variant, but, on the other hand, it results in a non-symmetric 
structure of the underlying elliptic sub-problems, which complicates the design efficient 
linear solver and preconditioning techniques. 

Instead of Problem~\ref{Prob:ContNewton} we then get following Newton iteration for the 
GGC$^1$(3) semidiscretization in time of the Navier--Stokes system 
\eqref{eq:navier_stokes} along with Nitsche's method for enforcing Dirichlet-type boundary 
conditions. Due to the cumbersome derivation of the equations and the innovation of the 
GGC$^1$(3) approach all formulas are explicitly given here in order to facilitate its 
application and implementation and enhance the confirmability of this work.

\begin{problem}[Newton iteration of GCC$^1$(3) time discretization with Nitsche's method]
	\label{Prob:ContNewton_N}
	Find corrections $(\delta \vec{v}_{n,2}^{k+1}, \delta \vec{v}_{n,3}^{k+1}) \in \mat{V}^2$ and $(\delta p_{n,2}^{k+1}, \delta p_{n,3}^{k+1})\in Q^2$ such that for all $\vec \phi = (\vec{\psi},\xi) \in \mat{V}^2 \times Q^2$ there holds that \begin{subequations}
	\begin{align*}
	\nonumber
	\mathcal{\widetilde{F}}_1 (\delta \vec{v}_{n,2}^{k+1}, \vec{\phi}) 
	&+ \mathcal{\widetilde{F}}_2 (\delta \vec{v}_{n,3}^{k+1}, \vec{\phi})
	+ \mathcal{\widetilde{B}}_1^{\,\top} (\delta p_{n,2}^{k+1}, \vec{\phi})
	+ \mathcal{\widetilde{B}}_2^{\,\top} (\delta p_{n,3}^{k+1}, \vec{\phi})
	\\[0.5ex]
	\nonumber
	& = - \Bigl<
	\vec{q}_1 (\vec{x}^k)
	\,, \vec{\psi} \Bigr>
	+
	\frac{\tau_n \nu}{2} \Bigl<
	\partial_{\vec{n}}\delta\vec{v}_{n, 2} (\vec{x}^k)
	\,, \vec{\psi} \Bigr>_{\Gamma_o}
	-
	\frac{\tau_n \nu}{12} \Bigl<
	\partial_{\vec{n}}\delta\vec{v}_{n, 3} (\vec{x}^k)
	\,, \vec{\psi} \Bigr>_{\Gamma_0}
	\\[0.5ex]
	%
	& \quad + \frac{\tau_n \nu}{2} \Bigl< \partial_{\vec{n}} \vec{\psi}, \vec{v}_{n, 2}^{k} - \vec{g}_{n, 2} \Bigr>_{\Gamma_D}
	- \frac{\eta_1 \tau_n}{2 h} \nu \Bigl< \vec{v}_{n, 2}^{k} - \vec{g}_{n, 2}, \vec{\psi} \Bigr>_{\Gamma_D}
	\\[0.5ex]
     \nonumber
	& - \frac{\eta_2 \tau_n}{2 h} \Bigl< (\vec{v}_{n, 2}^{k} - \vec{g}_{n, 2}) \cdot \vec{n}, \vec{\psi} \cdot \vec{n} \Bigr>_{\Gamma_D}
	- \frac{\tau_n \nu}{12} \Bigl< \partial_{\vec{n}} \vec{\psi}, \vec{v}_{n, 3}^{k} - \vec{g}_{n, 3} \Bigr>_{\partial \Gamma}
	\\[0.5ex]
	\nonumber
	& \quad + \frac{\eta_1 \tau_n}{12 h} \nu \Bigl< \vec{v}_{n, 3}^{k} - \vec{g}_{n, 3}, \vec{\psi} \Bigr>_{\partial \Gamma_D}
	+ \frac{\eta_2 \tau_n}{12 h} \Bigl< (\vec{v}_{n, 3}^{k} - \vec{g}_{n, 3}) \cdot \vec{n}, \vec{\psi} \cdot \vec{n} \Bigr>_{\Gamma_D} \,,
	\\[1.5ex]
	%
	- \mathcal{\widetilde{B}}_1 (\delta \vec{v}_{n,2}^{k+1}, \xi) 
	& - \mathcal{\widetilde{B}}_2 (\delta \vec{v}_{n,3}^{k+1}, \xi) 
	 =
	- \Bigl<
	q_2 (\vec{x}^k)
	\,, \xi \Bigr>
	- \frac{\tau_n}{2}\Bigl< \xi \vec{n}, \vec{v}_{n, 2}^{k} - \vec{g}_{n, 2} \Bigr>_{\Gamma_D}
	+ \frac{\tau_n}{12}\Bigl< \xi \vec{n}, \vec{v}_{n,3}^{k} - \vec{g}_{n, 3} \Bigr>_{\Gamma_D} \,,
	\\[1.5ex]
	\nonumber
	\mathcal{\widetilde{F}}_3 (\delta \vec{v}_{n,2}^{k+1},  \vec{\phi}) 
	& + \mathcal{F}_4 (\delta \vec{v}_{n,3}^{k+1},  \vec{\phi})
	+ \mathcal{\widetilde{B}}_3^{\,\top} (\delta p_{n,3}^{k+1},  \vec{\phi})
	=
	- \Bigl<
	\vec{q}_3 (\vec{x}^k)
	\,, \vec{\psi} \Bigr>
	\\[0.5ex]
	& - \nu \Bigl<
	\partial_{\vec{n}}\delta\vec{v}_{n, 2} (\vec{x}^k)
	\,, \vec{\psi} \Bigr>_{\Gamma_o}
	+ \nu \Bigl< \partial_{\vec{n}} \vec{\psi}, \vec{v}_{n, 2}^{k} - \vec{g}_{n, 2} \Bigr>_{\Gamma_D}
	- \frac{\eta_1}{h} \nu \Bigl< \vec{v}_{n, 2}^{k} - \vec{g}_{n, 2}, \vec{\psi} \Bigr>_{\Gamma_D}
	\\[0.5ex]
	\nonumber
	& - \frac{\eta_2}{h} \Bigl< (\vec{v}_{n, 2}^{k} - \vec{g}_{n, 2}) \cdot \vec{n}, \vec{\psi} \cdot \vec{n} \Bigr>_{\Gamma_D} \,,
	\\[1.5ex]
	%
	-\mathcal{\widetilde{B}}_3 (\delta \vec{v}_{n,3}^{k+1}, \xi) 
	& =
	- \Bigl<
	q_4 (\vec{x}^k)
	\,, \xi \Bigr>
	- \Bigl< \xi \vec{n}, \vec{v}_{n, 3}^{k} - \vec{g}_{n, 3} \Bigr>_{\Gamma_D}
	\,,
	\end{align*}	
    \end{subequations}
where the semi-linear and linear forms are defined by
\begin{align*}
&\begin{aligned}
\mathcal{\widetilde{F}}_1&(\delta \vec{v}_{n,2}^{k+1}, \vec{\phi}) \coloneqq
\Bigl< \delta \vec{v}_{n,2} ^{k+1} , \vec{\psi} \Bigr>
+ \frac{\tau_n \nu}{2} \Bigl<
\nabla \delta \vec{v}_{n,2}^{k+1}
, \nabla \vec{\psi} \Bigr>
+ \tau_n \Bigl<
\nabla \delta \vec{v}_{n,2}^{k+1} \vec{a}(\vec{x}^k)
, \vec{\psi} \Bigr> + \tau_n \Bigl<
\nabla \vec{a}(\vec{x}^k) \delta \vec{v}_{n,2}^{k+1}
, \vec{\psi} \Bigr>\\[0.7ex] 
&
- \frac{\tau_n \nu}{2} \Bigl< \partial_{\vec{n}} \delta \vec{v}_{n,2}^{k+1}, \vec{\psi} 
\Bigr>_{\Gamma_D}
- \frac{\tau_n \nu}{2} \Bigl< \partial_{\vec{n}} \vec{\psi}, \delta \vec{v}_{n, 2}^{k+1} 
\Bigr>_{\Gamma_D}
+ \frac{\eta_1 \tau_n}{2 h} \nu \Bigl< \delta \vec{v}_{n, 2}^{k+1}, \vec{\psi} 
\Bigr>_{\Gamma_D}
+ \frac{\eta_2 \tau_n}{2 h} \Bigl< \delta \vec{v}_{n, 2}^{k+1} \cdot \vec{n}, \vec{\psi} 
\cdot \vec{n} \Bigr>_{\Gamma_D} \,,
\end{aligned}
\end{align*}
\begin{align*}
&\begin{aligned}
\mathcal{\widetilde{F}}_2&(\delta \vec{v}_{n,3}^{k+1}, \vec{\phi}) \coloneqq
- \frac{\tau_n \nu}{12} \Bigl<
\nabla \delta \vec{v}_{n,3}^{k+1}
, \nabla \vec{\psi} \Bigr>
+ \tau_n \Bigl<
\nabla \delta \vec{v}_{n,3}^{k+1} \vec{b}(\vec{x}^k)
, \vec{\psi} \Bigr> + \tau_n \Bigl<
\nabla \vec{b}(\vec{x}^k) \delta \vec{v}_{n,3}^{k+1}
, \vec{\psi} \Bigr>
\\[0.7ex]
& + \frac{\tau_n \nu}{12} \Bigl< \partial_{\vec{n}} \delta \vec{v}_{n, 3}^{k+1}, \vec{\psi} \Bigr>_{\Gamma_D}
+ \frac{\tau_n \nu}{12} \Bigl< \partial_{\vec{n}} \vec{\psi}, \delta \vec{v}_{n, 3}^{k+1} \Bigr>_{\Gamma_D}  - \frac{\eta_1 \tau_n}{12 h} \nu \Bigl< \delta \vec{v}_{n, 3}^{k+1}, \vec{\psi} \Bigr>_{\Gamma_D}
- \frac{\eta_2 \tau_n}{12 h} \Bigl< \delta \vec{v}_{n, 3}^{k+1} \cdot \vec{n}, \vec{\psi} \cdot \vec{n} \Bigr>_{\Gamma_D} \,,
\end{aligned}
\\[1.5ex]
&\begin{aligned}
\mathcal{\widetilde{F}}_3&(\delta \vec{v}_{n,2}^{k+1}, \vec{\phi}) \coloneqq
\nu \Bigl<
\nabla \delta\vec{v}_{n,2}^{k+1}
, \nabla \vec{\psi} \Bigr>
+ \Bigl<
\nabla \vec{v}_{n,2}^k \delta\vec{v}_{n,2}^{k+1}
, \vec{\psi} \Bigr> + \Bigl<
\nabla \delta\vec{v}_{n,2}^{k+1} \vec{v}_{n,2}^k
, \vec{\psi} \Bigr>
- \nu \Bigl< \partial_{\vec{n}} \delta \vec{v}_{n, 3}^{k+1}, \vec{\psi} \Bigr>_{\Gamma_D}
\\[0.7ex]
& - \nu \Bigl< \partial_{\vec{n}} \vec{\psi}, \delta \vec{v}_{n, 2}^{k+1} \Bigr>_{\Gamma_D}+ \frac{\eta_1}{h} \nu \Bigl< \delta \vec{v}_{n, 2}^{k+1}, \vec{\psi} \Bigr>_{\Gamma_D}
+ \frac{\eta_2}{h} \Bigl< \delta \vec{v}_{n, 2}^{k+1} \cdot \vec{n}, \vec{\psi} \cdot \vec{n} \Bigr>_{\Gamma_D} \,,
\end{aligned}
\end{align*}
and, with $\mathcal{B}_i$ as well as $\mathcal{B}_i^{\,\top}$, $i = 1\ldots3$, as defined 
in \eqref{eq:bilinear_forms},
\begin{equation*}
\begin{aligned}
\mathcal{\widetilde{B}}_1 (\delta \vec{v}_{n,2}^{k+1}, \xi) & \coloneqq
\mathcal{B}_1
+ \frac{\tau_n}{2}\Bigl< \xi \vec{n}, \delta \vec{v}_{n, 2}^{k+1} \Bigr>_{\Gamma_D}, 
\quad
& \mbox{}\hspace*{-2ex}
\mathcal{\widetilde{B}}_1^\top (\delta p_{n,2}^{k+1}, \vec{\phi}) & \coloneqq
\mathcal{B}_1^\top
- \frac{\tau_n}{2} \Bigl< \delta p_{n,2}^{k+1} \vec{n}, \vec{\psi} \Bigr>_{\Gamma_D},
\\[1.5ex]
\mathcal{\widetilde{B}}_2 (\delta \vec{v}_{n,2}^{k+1}, \xi) &\coloneqq
\mathcal{B}_2
- \frac{\tau_n}{12}\Bigl< \xi \vec{n}, \delta \vec{v}_{n,3}^{k+1} \Bigr>_{\Gamma_D},
&  \mbox{}\hspace*{-2ex}
\mathcal{\widetilde{B}}_2^\top (\delta p_{n,3}^{k+1}, \vec{\phi})& \coloneqq
\mathcal{B}_2^\top
+ \frac{\tau_n}{12} \Bigl< \delta p_{n,3}^{k+1} \vec{n}, \vec{\psi} \Bigr>_{\Gamma_D} ,
\\[1.5ex]
\mathcal{\widetilde{B}}_3 (\delta \vec{v}_{n,2}^{k+1}, \xi) &\coloneqq
\mathcal{B}_3
+ \Bigl< \xi \vec{n}, \delta \vec{v}_{n, 3}^{k+1} \Bigr>_{\Gamma}\,, 
&  \mbox{}\hspace*{-2ex}
\mathcal{\widetilde{B}}_3^{\,\top} (\delta p_{n,3}^{k+1}, \vec{\phi}) & \coloneqq
\mathcal{B}_3^{\,\top}
- \Bigl< \delta p_{n,3}^{k+1} \vec{n}, \vec{\psi} \Bigr>_{\Gamma_D} \,.
\end{aligned}
\end{equation*}
\end{problem}

The fully discrete counterpart of Problem~\ref{Prob:ContNewton_N} is obtained along 
the lines of Subsec.~\ref{Subsec:SpaceDisc} with 
\begin{equation}
\label{eq:nitsche_matrices}
\begin{aligned}
\mat{\widetilde{F}}_1 & \coloneqq \left( \mathcal{\widetilde{F}}_1(\vec{\psi}_i, \vec{\psi}_j) \right)_{i,j=1}^J \,,
&
\mat{\widetilde{F}}_2 & \coloneqq \left( \mathcal{\widetilde{F}}_2(\vec{\psi}_i, \vec{\psi}_j) \right)_{i,j=1}^J \,,
&
\mat{\widetilde{F}}_3 & \coloneqq \left( \mathcal{\widetilde{F}}_3(\vec{\psi}_i, \vec{\psi}_j) \right)_{i,j=1}^J \,,
\\
\mat{\widetilde{B}}_1 & \coloneqq \left( \mathcal{\widetilde{B}}_1(\vec{\psi}_i, \vec{\psi}_j) \right)_{i, j  = 1}^J \,,
&
\mat{\widetilde{B}}_2 & \coloneqq \left( \mathcal{\widetilde{B}}_2(\vec{\psi}_i, \vec{\psi}_j) \right)_{i, j  = 1}^J \,,
&
\mat{\widetilde{B}}_3 & \coloneqq \left( \mathcal{\widetilde{B}}_3(\vec{\psi}_i, \vec{\psi}_j) \right)_{i, j  = 1}^J \,,
\\
\mat{\widetilde{B}}_1^\top & \coloneqq \left( \mathcal{\widetilde{B}}_1^\top(\vec{\psi}_i, \vec{\psi}_j) \right)_{i, j  = 1}^J \,,
&
\mat{\widetilde{B}}_2^\top & \coloneqq \left( \mathcal{\widetilde{B}}_2^\top(\vec{\psi}_i, \vec{\psi}_j) \right)_{i, j  = 1}^J \,,
&
\mat{\widetilde{B}}_3^\top & \coloneqq \left( \mathcal{\widetilde{B}}_3^\top(\vec{\psi}_i, 
\vec{\psi}_j) \right)_{i, j  = 1}^J 
\end{aligned}
\end{equation}
replacing the corresponding quantities \eqref{Def:F_iB_i}. The resulting block system of 
the Newton iteration has the same sparsity pattern as in 
\eqref{eq:newton_system_matrix}, but is now based on the matrices 
\eqref{eq:nitsche_matrices}.

\section{Numerical experiments}
\label{Sec:NumExp}

In this section we study numerically the GCC$^1$(3) approach along with Nitsche's 
method for Dirichlet boundary conditions, presented before in Sec.~\ref{Subsec:Nitsche}. 
Firstly, this is done by a numerical convergence study. A study of the condition number of the arising linear algebraic systems is also included. Secondly, the GCC$^1$(3) approach is applied to one of the popular benchmark problems proposed in \cite{TS96} of flow around a cylinder. The drag and lift coefficient are computed as goal quantities of physical interest. For the sake of comparion, calculations with the standard continuous Petrov--Galerkin method of piecewise linear polynomials in time, referred to as cGP(1), are also presented. For the implementation we used the deal.II finite element library \cite{arndt_daniel_deal.ii_2019} along with the Trilinos library  \cite{Trilinos-Users-Guide} for parallel computations on multiple processors.

\subsection{Convergence study}
\label{sec:convergence_study}

For the solution $\{\vec{v}, p\}$ of the Navier--Stokes system \eqref{eq:navier_stokes} 
and its fully discrete GCC$^1$(3) approximation $\{\vec{v}_{\tau,h}, p_{\tau,h}\}$ we 
define the errors 
\begin{align*}
\vec e^{\vec{v}} (t)&\coloneqq \vec{v}(t) - \vec{v}_{\tau,h}(t)\,,
&
e^{p}(t) &\coloneqq p(t) - p_{\tau,h}(,t)\,.
\end{align*}
We study the error $(\vec e^{\vec{v}},e^p)$ with respect to the norms
\begin{equation*}
\begin{aligned}
\| e^w \|_{L^\infty(L^2)} &\coloneqq  \max_{t \in I} \left( \int_{\Omega} \| e^w \|^2 \d 
x \right)^{\frac{1}{2}} \,,
\quad
\| e^w \|_{L^2(L^2)} \coloneqq  \biggl( \int_{I} \int_{\Omega} \| e^w(t) \|^2 \d x  \, \d t 
\biggr)^{\frac{1}{2}}\,,
\end{aligned}
\end{equation*}
where $w \in (\vec{v}, p)$.
The $L^{\infty}$-norm in time is computed on the discrete time grid
\begin{equation*}
I = \big\{
t_n^d \mid t_n^d = t_{n-1} + d \cdot k_n \cdot \tau_n,
\quad
k_n = 0.001\,,\; d = 0, \ldots, 999\,,\; n = 1, \ldots, N
\big\}\,.
\end{equation*}

In our experiment we study a test setting presented in \cite{BR12} and choose 
the right-hand side function $\vec{f}$ on $\Omega \times I = (0,1)^2 \times (0,1]$ in 
such a way, that the exact solution  of the Navier--Stokes system~\eqref{eq:navier_stokes} 
is given by
\begin{equation}
\label{Eq:AnalSol}
\begin{aligned}
\vec{v}(\vec{x},t)
&\coloneqq
\begin{pmatrix}
\cos(x_2 \pi) \cdot \sin(t) \cdot \sin(x_1 \pi)^2  \cdot \sin(x_2 \pi) \\[0.5ex]
- \cos(x_1 \pi) \cdot \sin(t) \cdot \sin(x_2 \pi)^2  \cdot \sin(x_1 \pi)
\end{pmatrix}\,,
\\[1ex]
p(\vec{x},t)
&\coloneqq
\cos(x_2 \pi) \cdot \sin(t) \cdot \sin(x_1 \pi) \cdot \cos(x_1 \pi) \cdot \sin(x_2 \pi)\,.
\end{aligned}
\end{equation}
We prescribe a Dirichlet boundary condition \eqref{eq:navier_stokes_3}, given by 
the solution \eqref{Eq:AnalSol}, on the whole boundary such that $\Gamma_D=\partial 
\Omega$, i.e.\ $\vec g = \vec 0$. The initial condition  \eqref{eq:navier_stokes_3} is 
also given by \eqref{Eq:AnalSol}, i.e.\ $\vec v_0 = \vec 0$. For the discretization in 
space the $\mathcal Q_4$--$\mathcal Q_3$ pair of the Taylor--Hood family is used; cf.\ 
Fig.~\ref{fig:sparsity_pattern}. The Nitsche penalty parameters in \eqref{eq:Nitsche_3} 
are fixed to $\eta_1 = \eta_2 = 35$.
We also compute the spectral condition number $\kappa_2(\mat{S})$ of the corresponding system matrices $\mat{S}$ of the Newton linearization, using the largest singular value $\sigma_1$ and its smallest one $\sigma_n$, by means of 
\begin{equation*}
    \kappa_2(\mat{S}) = \frac{\sigma_1}{\sigma_n} \,.
\end{equation*}
\sisetup{scientific-notation = true,
		 round-mode=places,
		 round-precision=3,
		 output-exponent-marker=\ensuremath{\mathrm{e}},
		 table-figures-integer=1, 
		 table-figures-decimal=3, 
		 table-figures-exponent=1, 
		 table-sign-mantissa = false, 
		 table-sign-exponent = false, 
	 	 table-number-alignment=center} 
\begin{table}[t]
    \caption{Step sizes, resulting degrees of freedom (DoF) on each time interval $I_n$, spectral condition number $\kappa_2(\mat{S})$, errors and experimental orders of convergence (EOC) for the approximation of \eqref{Eq:AnalSol} with $\mathcal Q_4$--$\mathcal Q_3$ elements in space and different time integration schemes under Nitsche's method; $\tau_0 = \num{1.0}$, $h_0 = 1/\sqrt 2$.}
	\centering
    \begin{subtable}{1.\textwidth}
    \centering
    \caption{$GCC^1(3)$ time integration scheme}
    \begin{tabular}{c@{\,\,\,\,}c  rS  S@{\,}S  S@{\,}S}
		\toprule
        {$\tau$} & {$h$} &
        {$\text{DoF}_{|I_n}$} & {$\kappa_2(\mat{S})$} & { $\| e^{\;\vec{v}}  \|_{L^2(L^2)} $ } &
        { $\| e^{\;p}  \|_{L^2(L^2)} $ } &
        { $\| e^{\;\vec{v}}  \|_{L^\infty(L^2)} $ } &
		{ $\| e^{\;p}  \|_{L^\infty(L^2)} $ } \\
        \cmidrule(r){1-2}
        \cmidrule(lr){3-4}
        \cmidrule(l){5-8}
        $\tau_0/2^0$  &  $h_0/2^0$  &  $422$      &  1.0845e6   &  3.099e-04  &  9.240e-04  &  7.082e-03  &  3.139e-03  \\
        $\tau_0/2^1$  &  $h_0/2^1$  &  $1\,494$   &  5.71573e6  &  1.954e-05  &  5.516e-05  &  4.437e-05  &  1.973e-04  \\
        $\tau_0/2^2$  &  $h_0/2^2$  &  $5\,606$   &  4.36986e7  &  1.226e-06  &  3.468e-06  &  2,780e-06  &  1.264e-05  \\
        $\tau_0/2^3$  &  $h_0/2^3$  &  $21\,702$  &  3.45893e8  &  7.673e-08  &  2.177e-07  &  1.734e-07  &  8.016e-07  \\
        \cmidrule(r){1-2}
        \cmidrule(lr){3-4}
        \cmidrule(l){5-8}
        \multicolumn{2}{c}{EOC} &
        \multicolumn{1}{c}{--} & \multicolumn{1}{c}{--} & {4.00} & {3.99} &
        {4.00} & {3.98} \\
		\bottomrule
        \addlinespace[2ex]
	\end{tabular}
    \label{tab:conv_gcc}
\end{subtable}%
\\%
\begin{subtable}{1.\textwidth}
    \centering
    \caption{cGP(1) time integration scheme}
    \begin{tabular}{c@{\,\,\,\,}c  rS  S@{\,}S  S@{\,}S}
        \toprule
        {$\tau$} & {$h$} &
        {$\text{DoF}_{|I_n}$} & {$\kappa_2(\mat{S})$} & { $\| e^{\;\vec{v}}  \|_{L^2(L^2)} $ } &
        { $\| e^{\;p}  \|_{L^2(L^2)} $ } &
        { $\| e^{\;\vec{v}}  \|_{L^\infty(L^2)} $ } &
		{ $\| e^{\;p}  \|_{L^\infty(L^2)} $ } \\
        \cmidrule(r){1-2}
        \cmidrule(lr){3-4}
        \cmidrule(l){5-8}
        $\tau_0/2^0$  &  $h_0/2^0$  &  $211$      &  7.21336e4  &  1.3583302742e-02  &  1.1702034258e-02  &  1.7899395044e-02  &  1.7899395044e-02  \\
        $\tau_0/2^1$  &  $h_0/2^1$  &  $747$      &  2.23797e5  &  3.6515714874e-03  &  3.3942495442e-03  &  5.6622334073e-03  &  5.6622334073e-03  \\
        $\tau_0/2^2$  &  $h_0/2^2$  &  $2\,803$   &  8.9551e5   &  9.2699994632e-04  &  8.2779775257e-04  &  1.5017817052e-03  &  1.5017817052e-03  \\
        $\tau_0/2^3$  &  $h_0/2^3$  &  $10\,851$  &  3.58716e6  &  2.3257321921e-04  &  2.0574424224e-04  &  3.8051175754e-04  &  3.8051175754e-04  \\
        \cmidrule(r){1-2}
        \cmidrule(lr){3-4}
        \cmidrule(l){5-8}
        \multicolumn{2}{c}{EOC} &
        \multicolumn{1}{c}{--} & \multicolumn{1}{c}{--} & {1.99} & {2.01} &
        {1.93} & {1.98} \\
        \bottomrule
    \end{tabular}
    \label{tab:conv_cgp}
\end{subtable}%
\label{tab:numerical_results_1}
\end{table}
\Cref{tab:conv_gcc} shows the calculated errors as well as the 
experimental orders of convergence and the spectral condition numbers for a sequence of meshes that are 
successively refined in space and time. We note that a test of 
simultaneous convergence in space and time is thus performed. In all measured norms, we 
observe convergence of fourth order. This is the optimal order for the 
Galerkin--collocation approach GCC$^1$(3) with piecewise polynomials of order three in 
time. For the mixed approximation of the Navier--Stokes system by the $\mathcal 
Q_4$--$\mathcal Q_3$ pair of the Taylor--Hood family convergence of order four in space 
can at most be expected. Thus, the application of the Nitsche method does not 
deteriorate the convergence behavior. The optimal rate of convergence in time is thus obtained for the approximation of the velocity field and of the pressure variable. 
For comparison and to illustrate the quality of the GCC$^1$(3) approach, Table~\ref{tab:conv_cgp} also shows the results for the standard cGP$(1)$ discretization in time of second order accuracy. Of course, comparing the step sizes or number of the degrees of freedom with the errors for both approximation schemes, the higher order GCC$^1$(3) method is superior to the cGP(1) one. Comparing the condition of the linear algebraic systems with the computed errors for both approximation schemes, we observe that for fixed condition numbers the higher order GCC$^1$(3) method yields smaller errors than the cGP(1) one. Consequently, the conditioning of the linear systems, as a measure for the complexity of their iterative solution, does not suffer from the higher order combined Galerkin--collocation approximation.

Finally, we note that time-dependent boundary conditions can be captured by the 
Galerkin--collocation approach without loss of order of convergence. This is illustrated  
numerically in \cite{AB19_1} for the wave equation. 

\subsection{Impact of Nitsche's method for flow around a cylinder}
\label{sec:flow_cyl}

\begin{figure}[htb!]
	\centering
	\includegraphics[width=0.7\textwidth,keepaspectratio]
    {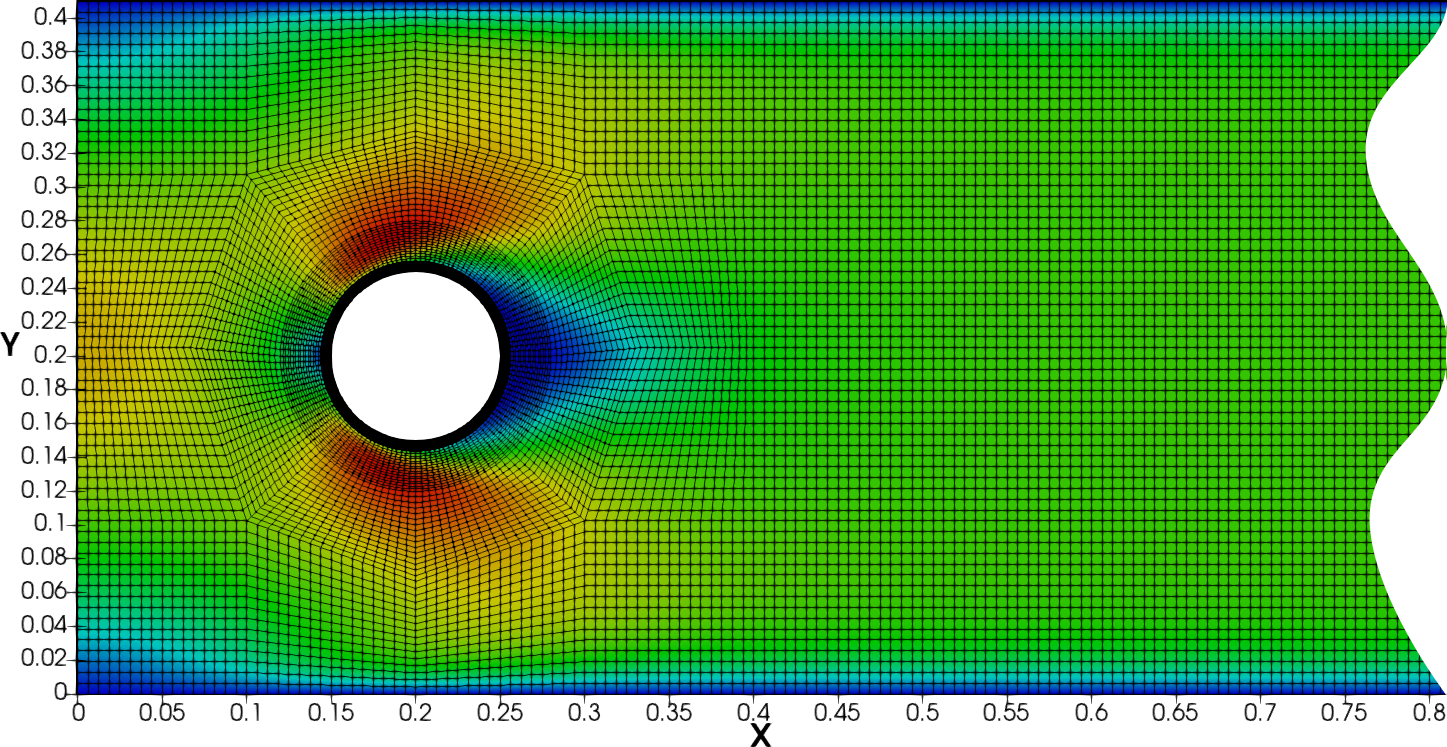}
	\caption{Inflow boundary at $x = 0$ and part of the spatial grid around the cylinder.}
	\label{fig:dfg_grid}
\end{figure}

In the second numerical example, we compare the effect of imposing the boundary 
conditions in a weak form by using Nitsche's method (cf.\ Subsec.~\ref{Subsec:Nitsche}) 
with enforcing the boundary conditions by the definition of the underlying function space 
(cf.\ Subsec.~\ref{Subsec:Semdisc} and \ref{Subsec:SpaceDisc}) and, then, condensing the 
algebraic system by eliminating the degrees of freedom corresponding to the nodes on the 
Dirichlet part of the boundary. For the experimental setting we use the well-known DFG 
benchmark problem "flow around a cylinder", given in \cite{TS96}. A section of the mesh used for the computations is shown in \cref{fig:dfg_grid}.

We consider the time interval $I = (0, 1]$, set $\nu = 0.01$ and let the velocity on the inflow boundary 
$\Gamma_i$ be given by ($\vec x = (x,y)^\top$)
\begin{align*}
\vec{g}(x,y,t)
=
\begin{pmatrix}
-7.13861\cdot (y - 0.41) \cdot y \cdot t^2 \\
0
\end{pmatrix} \,.
\end{align*}
The maximum mean velocity of the parabolic inflow profile is reached for $T=1$ and 
is $\bar{U} = 0.2$. With the diameter of the cylinder as the characteristic length $L = 
0.1$, this results in a Reynolds number $Re$ of
\begin{equation*}
	Re = \frac{\bar{U} \cdot L}{\nu} = \frac{0.2 \cdot 0.1}{0.01} = 2 \,.
\end{equation*}
\begin{figure}[!ht]
	\centering
	\includegraphics[width=0.48\textwidth]
	{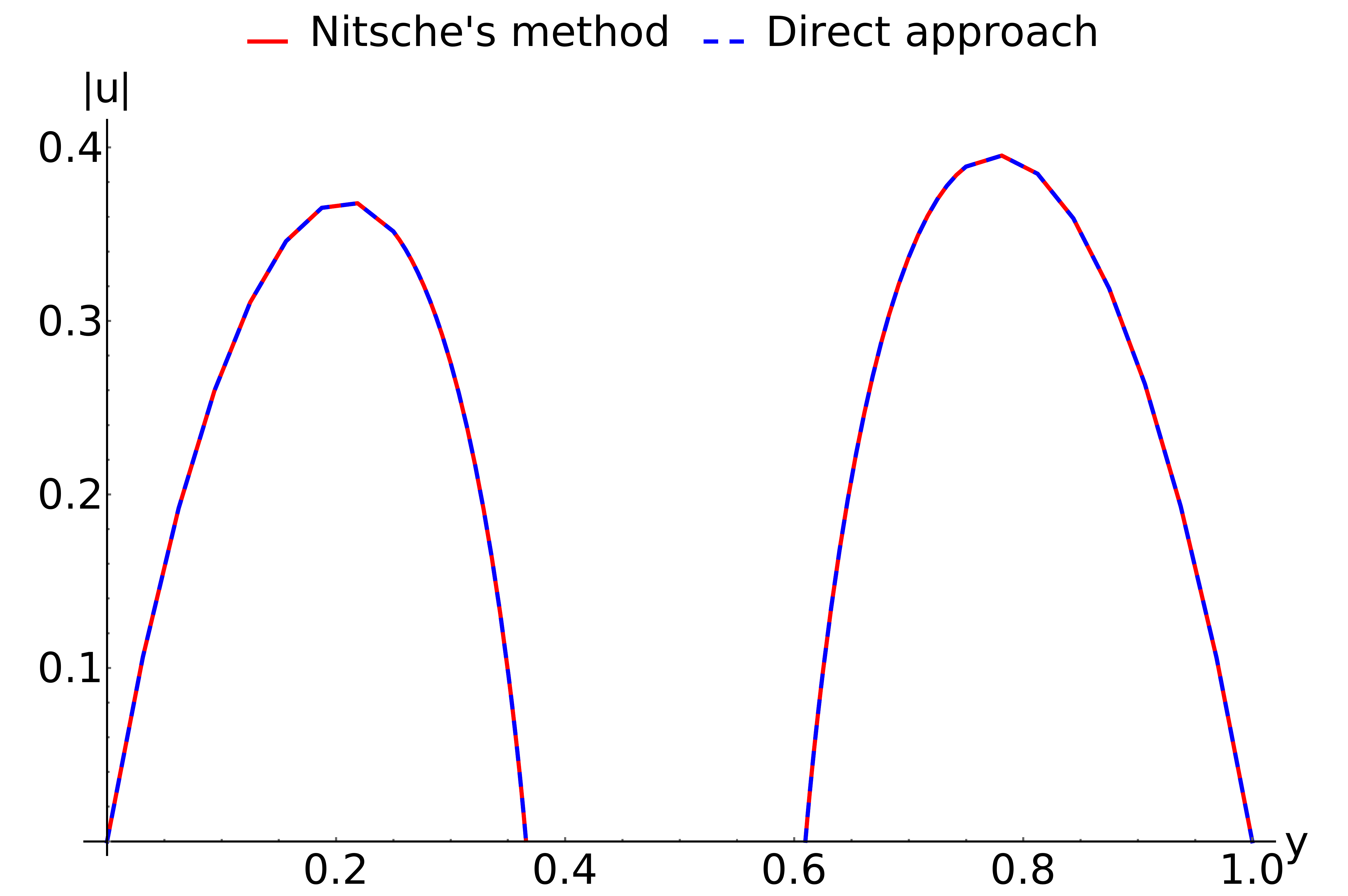}
	\hspace*{2ex}
	\includegraphics[width=0.48\textwidth]
	{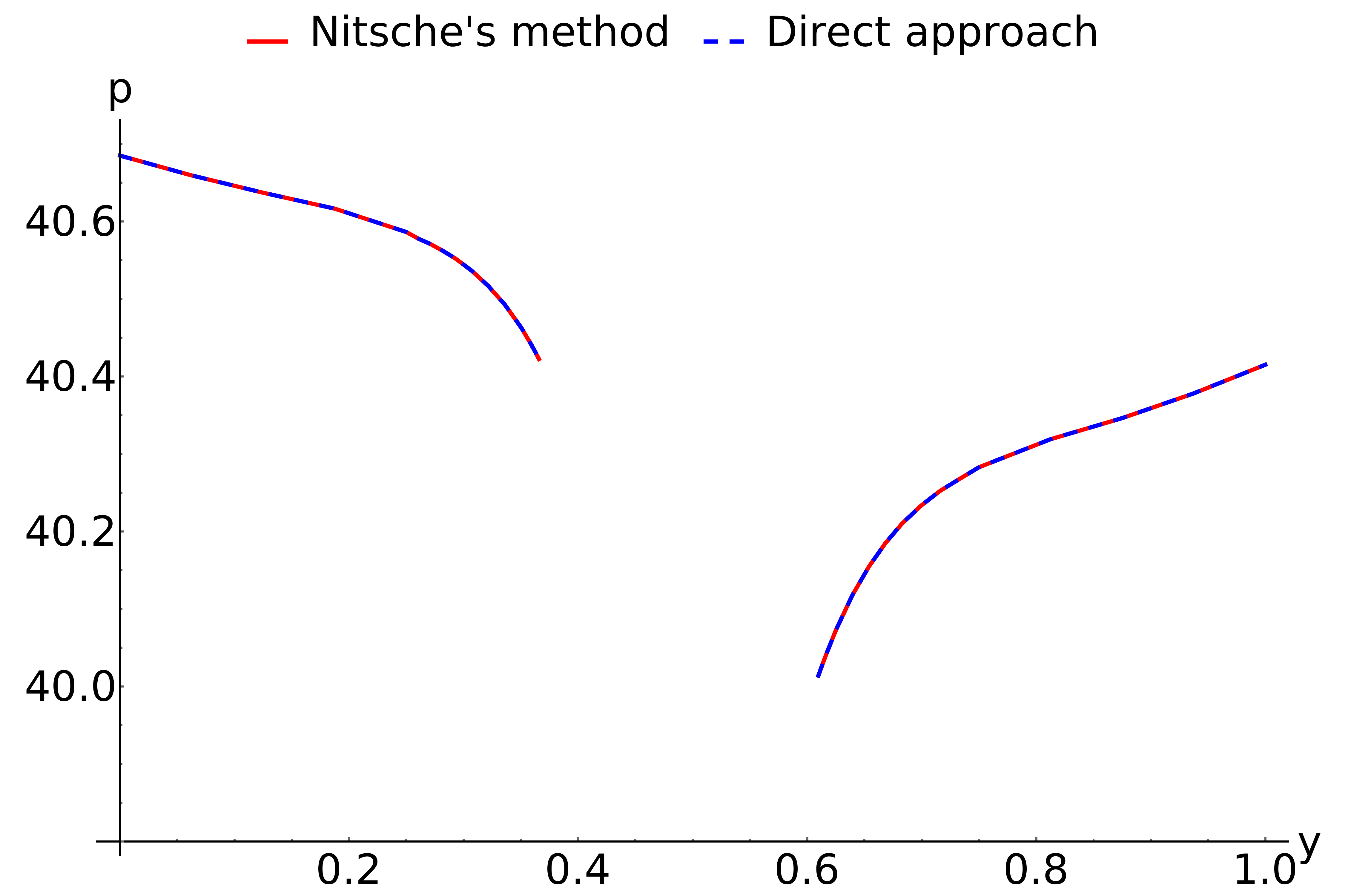}
	\caption{Comparison of computed velocity and pressure profiles along the 
$y$-axis at $x=0.2$ (cross section line through the ball's midpoint) and $t = 1$ for 
Nitsche's method (weak form of Dirichlet boundary conditions) and for the enforcement of 
the Dirichlet boundary conditions by the definition of the function spaces with 
condensation of the algebraic system.}
	\label{fig:v_p_comparison}
\end{figure}
In Fig.\ \ref{fig:v_p_comparison} we compare the computed solutions along the 
$y$-direction at $x = 0.2$ (cross section line through the ball's midpoint) and $T = 1$, 
that are obtained by the either methods. In the figures, the interval without 
having graphs is the cross section line that is covered by the ball. The computed 
profiles match perfectly such that no loss of accuracy is observed by the application of 
Nitsche's method of enforcing Dirichlet boundary conditions for this problem of viscous 
flow.

\subsection{GCC$^1(3)$ versus cGP(1) for time-periodic flow around a cylinder}
\label{sec:flow_re_100}


In the third numerical example we evaluate the performance properties of the Galerkin--collocation approach GCC$^1(3)$ for the DFG benchmark (cf.\ \cite{TS96}) of the time-periodic behavior of a fluid in a pipe with a circular obstacle as a more sophisticated test problem. It is set up in two space dimensions. The drag and lift coefficients of the flow on the circular cross section are computed as goal quantities of physical interest. We compare the computed results of the GCC$^1(3)$ scheme with the ones obtained by the cGP(1) approximation (or Crank--Nicholson method). We use the geometrical setup as precribed in \cref{sec:flow_cyl}, but consider the time interval $I = (0, 10]$, set the initial time step size $\tau = 0.01$, the viscosity $\nu = 0.001$ and prescribe the following Dirichlet boundary condition on $\Gamma_i$ by
\begin{equation}
\label{Eq:DFGBM_In}
    g(x, y, t) =
    \begin{pmatrix}
        \frac{4 \cdot 1.5 \cdot y (0.41 - y)}{0.41^2} 
        \cdot \Bigl( (3 t^2 - 2 t^3) \cdot \bigl( 1 - \Theta(t - 1) \bigr)
        + \Theta(t-1) \Bigr)
        \\[1ex]
        0
    \end{pmatrix} \,,
\end{equation}
where $\Theta$ is the Heaviside function. Only for reasons of convenience and implementational issues, we altered the time step size slightly by choosing it as a constant, compared with the original benchmark design in \cite{TS96}.
%
%
Eq.\ \eqref{Eq:DFGBM_In} coincides with the inflow condition given in \cite[p.\,4]{TS96}, except that a smooth and non-instantaneous increase in time of the profile until $t=1$ is used. For $t \geq 1$ this results in a Reynold's number of $Re = 100$ and a time-periodic flow behaviour. With the drag and lift forces $F_D$ and $F_L$ on the circle $S$ that are defined by
\begin{align}
    F_D &= \int_S
    \left(
        \nu \frac{\partial v_t}{\partial \vec{n}} n_y - P n_x
    \right) \d S \,,
        &
    F_L &= - \int_S
    \left(
        \nu \frac{\partial v_t}{\partial \vec{n}} n_x - P n_y
    \right) \d S \,,
\end{align}
where $\vec n$ is the normal vector on $S$ and $v_t$ is the tangential velocity along $\vec{t} = (n_y, -n_x)^\top$, we can compute the drag and lift coefficients $c_D, c_L$ as our goal quantities by
\begin{align}
    c_D &= \frac{2}{\bar{U}^2 L} F_D \,,
        &
    c_L &= \frac{2}{\bar{U}^2 L} F_L \,.
\end{align}
\begin{figure}[H]
	\centering
	\subcaptionbox{Drag coefficients $C_D$ \label{fig:drag_coefficients}}
	[1.\columnwidth]
	{\includegraphics[width=0.82\textwidth,keepaspectratio]{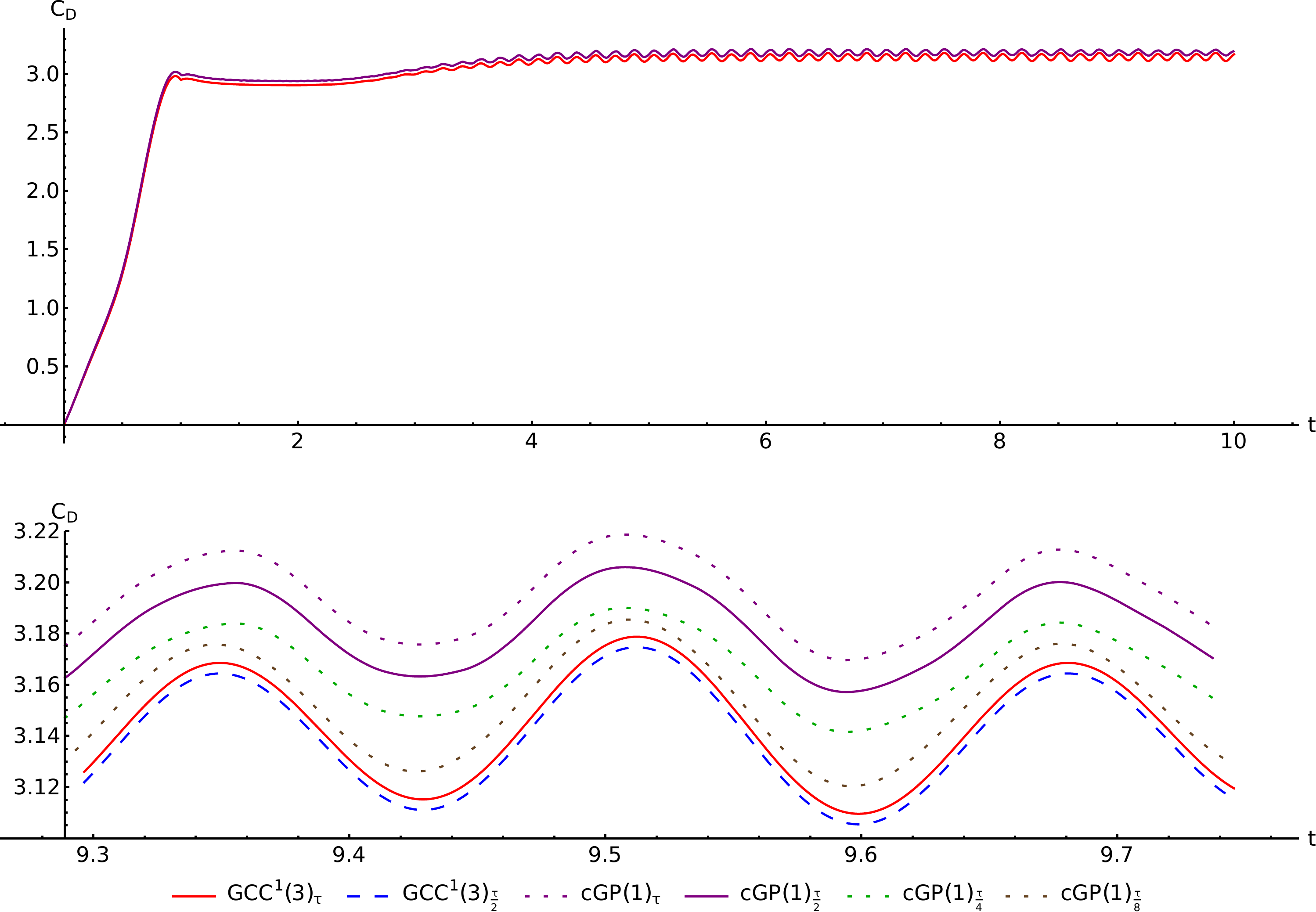}
	}\\[5ex]
	\subcaptionbox{Lift coefficients $c_L$ \label{fig:lift_coefficients}}
	[1.\columnwidth]
	{\includegraphics[width=0.82\textwidth,keepaspectratio]{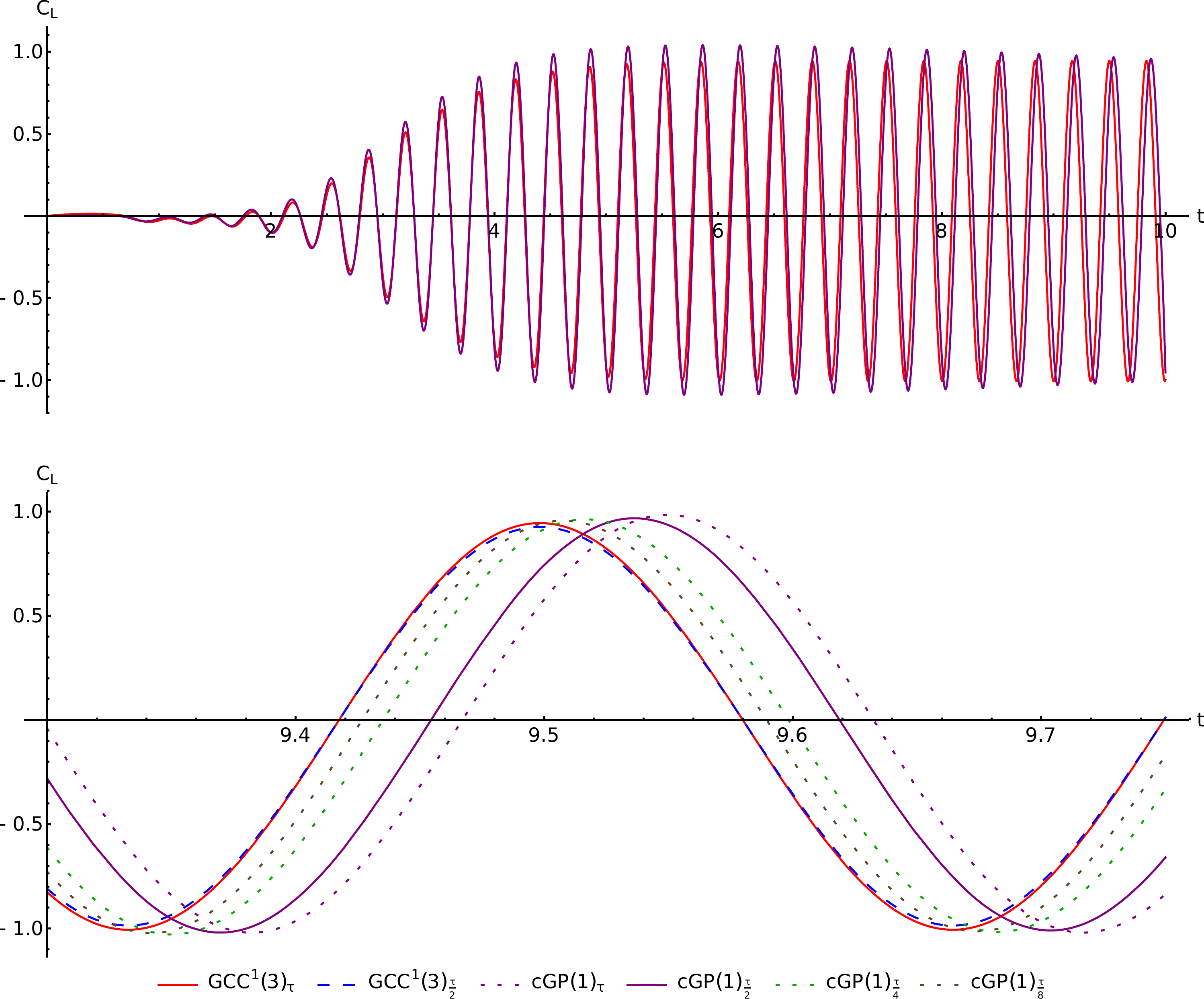} 
	}
	\caption{Computed drag and lift coefficients for the example of \cref{sec:flow_re_100} and different time discretization schemes with basic time step size $\tau$.}
	\label{fig:drag_lift}
\end{figure}

For the spatial discretization we use a mesh consisting of $6\,912$ cells with higher order $\mathcal Q_5$--$\mathcal Q_4$ Taylor--Hood elements. This results in $919\,104$ degrees of freedom in each time interval for the GCC$^1(3)$ scheme in contrast to $459\,552$ degrees of freedom for the cGP$(1)$ method. In Fig.~\ref{fig:drag_lift} we visualize the computed values for the drag and lift coefficient. The solid lines represent GCC$^1(3)$ and cGP$(1)$ simulations with same number of degrees of freedom, summarized over the whole simulation time. It can be observed that the cGP$(1)$ solution converges towards the GCC$^1(3)$ solution. But even with an eighth of the time step size of the GCC$^1(3)$ solution it does not completely coincide with the GCC$^1(3)$ solution yet. In contrast to this, the GCC$^1(3)$ solution almost seems to be fully converged in the time domain with the basic time step size of $\tau$, since its further reduction to $\frac{\tau}{2}$ results in a minor change of the resulting drag and lift coefficients only. Summarizing, we can state that the GCC$^1(3)$ approach is strongly superior to the cGP$(1)$ method with respect to accuracy and efficiency.

\section{Outlook}
\label{Sec:Outlook}

We proposed and analyzed numerically higher order Galerkin--collocation time discretization schemes along with Nitsche's method for incompressible viscous flow. The time discretization combines Galerkin approximation with the concepts of collocation. Expected optimal order convergence properties were obatined in a numerical experiment. A careful comparative numerical study with the continuous Galerin--Petrov method with piecewise linear polynomials was presented further for the DFG benchmark of flow around a cylinder with $Re=100$. The higher order Galerkin--collocation approach offered the potential of usage of much larger time steps withput loss of accuracy and, thus, improves the efficiency of the time discretization strongly. The solver of the resulting linear systems (cf.\ Eq.~\eqref{eq:newton_system_matrix}) still continues to remain an important research topic for the future. A geometric 
multigrid preconditioner, based on a Vanka smoother, for GMRES outer iterations, showed promising results for the Navier--Stokes equations; cf. \cite{HST14}. A similar approach was successfully used for simulations in three space dimensions of fully coupled fluid-structure interaction problems \cite{FR19}. It remains to improve the current 
solver and develop a similar competitive solver and preconditioner for the Galerkin--collocation approach to the Navier--Stokes system. This will be a work for the future.

\begin{figure}[h!tb]
	\centering
	\includegraphics[width=0.7\textwidth]
		{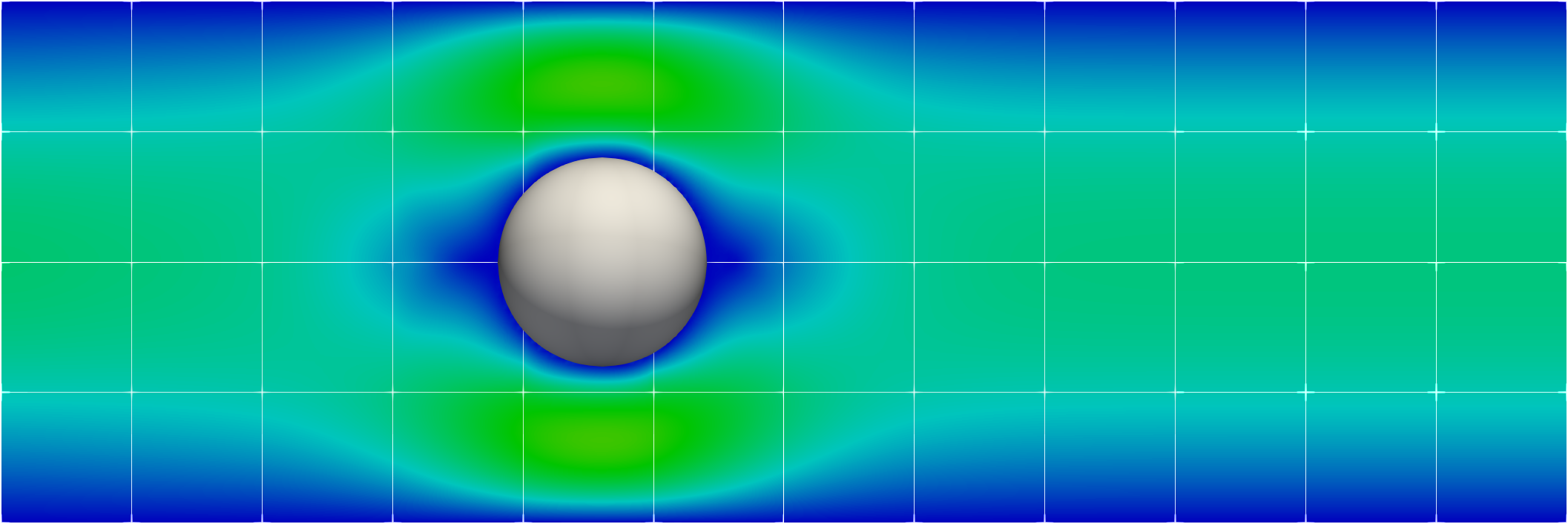}\\[2ex]
	\includegraphics[width=0.7\textwidth]
		{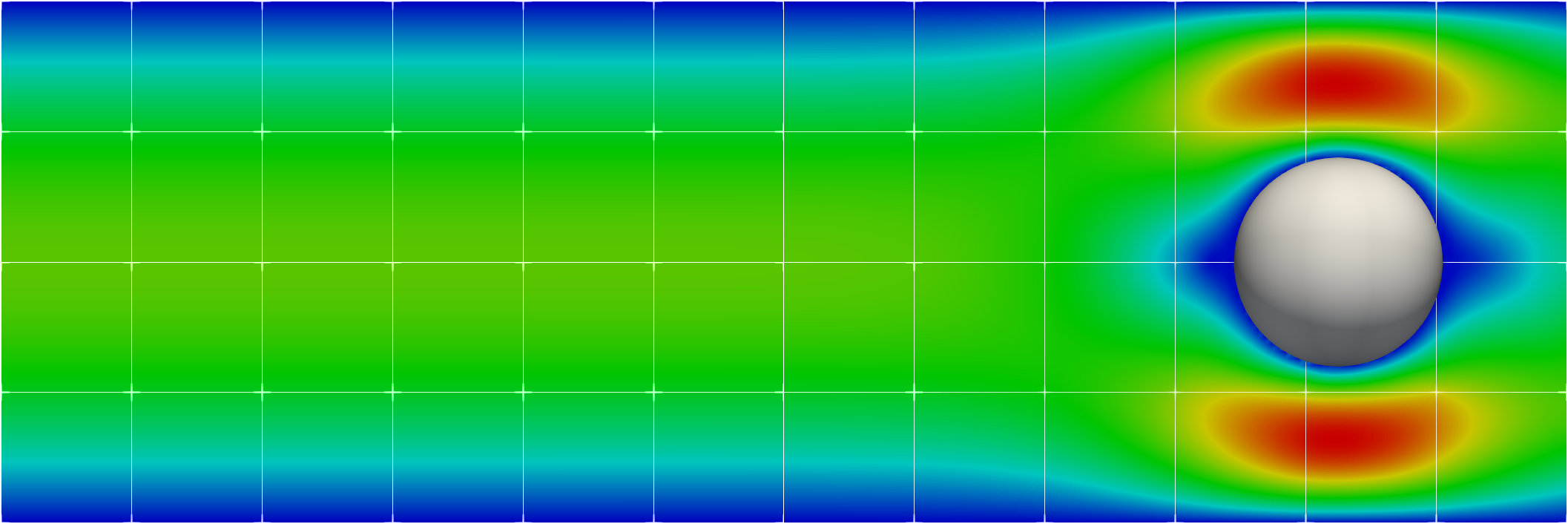}
	\caption{Application of Nitsche's method with cut finite elements on a dynamic 
geometry.}
		\label{Fig:MovBall}
\end{figure}

Enforcing  Dirichlet boundary conditions in a variational form by Nitsche's method offers the potential of applying the Galerkin--collocation approach to complex flow problems in dynamic geometries with moving boundaries or interfaces. In particular, fluid structure interaction is in the scope of our interest by combining the concepts of this work with our former work \cite{AB19_1,AB19_2} on the numerical approximation of wave phenomena. The capability of treating dynamic domains still requires the application of the concept of fictitious domains, that is based on Nitsche's method, together with using cut finite element techniques; cf.\ \cite{BH14,MLLR14,MSW18,S17} and the references therein. As a proof of concept, we illustrate in Fig.~\ref{Fig:MovBall} a preliminary result for the simulation of flow around a moving ball. This demonstrates that the proposed Galerkin--collocation approach along with 
Nitsche's method offers high potential for problems on evolving domains. In Fig.~\ref{Fig:MovBall}, the flow around a ball, that is moving forward and backward with prescribed velocity in a pipe, has been computed by using the Galerkin--collocation approach GCC$^1$(3) along with a Nitsche fictitious domain method and cut finite elements; cf.\ \cite{AB20}. The background mesh is kept fixed for the whole simulation time such that cut elements arise, Further details will be presented in a forthcoming work since it would overburden this paper.




\end{document}